\documentclass[twoside,reqno]{amsart} 

\usepackage{latexsym,rotate,eucal,cite,url}
\usepackage{amsmath,amsthm,amssymb,amsxtra} 

\newcommand{\quab}{\hspace*{2.0mm}}
\newcommand{\puab}{\hspace*{-1.5mm}}
\newcommand{\xpuad}{\hspace*{-9.0mm}}

\newcommand{\qqed}{\tag*{\qed}} 

\newcommand{\pav}[1]{\lfloor{#1}\rfloor}

\newcommand{\pavv}[1]{\left\lfloor{#1}\right\rfloor}

\newcommand{\ang}[1]{\left\langle{#1}\right\rangle}

\newcommand{\mb}[1]{\mathbb{#1}}
\newcommand{\mc}[1]{\mathcal{#1}}
\newcommand{\md}{\mathrm}

\newcommand{\ff}[2]{(#1)_{#2}}
\newcommand{\hyp}[3]{\left[\puab\ba{#1}#2\\#3\ea\puab\right]}
\newcommand{\hyq}[4]{\left[\puab\ba{#1}#3\\#4\ea{\!\!\Big|#2}\right]}

\newcommand{\binm}{\binom}
\newcommand{\binq}[2]{{#1\brack#2}}

\newcommand{\be}{\begin{equation}}
\newcommand{\ee}{\end{equation}}
\newcommand{\ba}{\begin{array}}
\newcommand{\ea}{\end{array}}
\newcommand{\bmn}{\begin{eqnarray}}
\newcommand{\emn}{\end{eqnarray}}
\newcommand{\bnm}{\begin{eqnarray*}}
\newcommand{\enm}{\end{eqnarray*}}
\newcommand{\bln}{\begin{subequations}}
\newcommand{\eln}{\end{subequations}}

\newcommand{\pq}[1]{\begin{equation}#1\end{equation}}
\newcommand{\pp}[2]{\begin{aligned}[#1]#2 
            \end{aligned}}  
\newcommand{\pmq}[1]{\begin{align}#1
            \end{align}}     
\newcommand{\pnq}[1]{\begin{align*}#1
            \end{align*}}    
\newcommand{\pnp}[2]{\begin{alignat*}{#1}#2
            \end{alignat*}}  
\newcommand{\mult}[2]{\begin{array}{#1}#2\end{array}}%

\newcommand{\centro}[1]
           {\begin{center}#1\end{center}}

\newcommand{\alp}{\alpha}
\newcommand{\bet}{\beta}
\newcommand{\gam}{\gamma}

\newcommand{\lam}{\lambda}

\newcommand{\veps}{\varepsilon}
\newcommand{\del}{\delta}

\newcommand{\sig}{\sigma}
\newcommand{\vph}{\varphi}

\newcommand{\Lam}{\Lambda}
\newcommand{\Ome}{\Omega}

\newcommand{\Gam}{\Gamma}

\newtheorem{thm}{Theorem}
\newtheorem{lemm}[thm]{Lemma}

\newtheorem{prop}[thm]{Proposition}
\newtheorem{exam}{Example}

\newtheorem{entry}{Entry}

\newcommand{\referxy}[4]{\bibitem{kn:#1}{#2,}~\emph{#3,}~{#4.}}	
\newcommand{\cito}[1]{\cite{kn:#1}}	
\newcommand{\citu}[2]{\cite[#2]{kn:#1}}

\begin{document}{\fbox{June 22, 2021}\hfill\fbox{\bf Ninth version}}%
\title{$q$-Analogues of $\pi$-Related Formulae from\\
       Jackson's $_8\phi_7$-Series via Inversion Approach}
\author{Xiaojing Chen and Wenchang Chu}
\address{School of Statistics\newline
	Qufu Normal University\newline
	Qufu (Shandong), P.~R.~China}
\email{upcxjchen@outlook.com}
\address{Department of Mathematics and Physics\newline
     	University of Salento (P.~O.~Box~193) \newline
      	73100 Lecce, ~ Italy}
\email{chu.wenchang@unisalento.it}
\thanks{Email addresses: chu.wenchang@unisalento.it
	and upcxjchen@outlook.com}
\subjclass{Primary 33D15, Secondary 33C20, 65B10}
\keywords{Basic hypergeometric series; Bisection series;
          Jackson's formula for well--poised $_8\phi_7$-series;
	  Duplicate inversions; Triplicate inversions;
          Ramanujan--like $\pi$-related series;
	  Guillera's infinite series for $1/\pi^2$}


\begin{abstract}
By making use of the multiplicate form of the extended Carlitz inverse
series relations, we establish two general `dual' theorems of Jackson's
summation formula for well--poised $_8\phi_7$-series. Their duplicate
forms under the partition pattern $n=\lfloor{\frac{n}2}\rfloor+\lfloor{\frac{n+1}2}\rfloor$
are explored and yield numerous $q$-series identities whose limiting
cases as $q\to1$ result in classical $\pi$-related Ramanujan--like
series of convergence rate ``$\frac1{16}$" including one for $1/\pi^2$
discovered by Guillera (2003). The triplicate dual formulae under the
partition pattern $n=\lfloor{\frac{n}3}\rfloor+\lfloor{\frac{n+1}3}\rfloor+\lfloor{\frac{n+2}3}\rfloor$
are examined via the ``reverse bisection method", which leads us
to twenty new $q$-series identities together with their classical
counterparts of convergence rate ``$\frac{-1}{27}$" when $q\to1$.
\end{abstract}

\maketitle\thispagestyle{empty}
\markboth{Xiaojing Chen and Wenchang Chu}{$q$-Analogues
	  of $\pi$-Related Formulae via Inversion Approach}

\section{Introduction and Motivation}
Let $\mb{N}$ and $\mb{N}_0$ be the sets of natural numbers
and nonnegative integers, respectively. For an indeterminate
$x$, the shifted factorial is defined by
\[(x)_0\equiv1\quad\text{and}\quad
(x)_n=x(x+1)\cdots(x+n-1)
\quad\text{for}\quad n\in\mb{N}\]
with the following shortened multiparameter notation
\[\hyp{cccc}{\alp,\bet,\cdots,\gam}{A,B,\cdots,C}_n
=\frac{(\alp)_n(\bet)_n\cdots(\gam)_n}{(A)_n(B)_n\cdots(C)_n}.\]
Analogously, the shifted factorials with the
base $q$ are given by $(x;q)_0\equiv1$ and
\[(x;q)_n=(1-x)(1-qx)\cdots(1-q^{n-1}x)
\quad\text{for}\quad n\in\mb{N}.\]
The product and quotient of the $q$-shifted
factorials will be abbreviated respectively to
\pnq{\ff{\alp,\bet,\cdots,\gam;q}{n}\quab
&=\ff{\alp;q}{n}\ff{\bet;q}{n}\cdots\ff{\gam;q}{n},\\
\hyq{cccc}{q}{\alp,\bet,\cdots,\gam}
          {A,B,\cdots,C}_n
&=\frac{\ff{\alp;q}{n}\ff{\bet;q}{n}\cdots\ff{\gam;q}{n}}
     {\ff{A;q}{n}\ff{B;q}{n}\cdots\ff{C;q}{n}}.}
When $|q|<1$, the infinite product $(x;q)_{\infty}$ is well defined,
from which we introduce the $q$-gamma function~\citu{rahman}{\S1.10}
\[\Gam_q(x)\:=\:(1-q)^{1-x}
\frac{(q;q)_\infty}{(q^x;q)_\infty}
\quad\text{and}\quad
\lim_{q\to1^{-}}\Gam_q(x)\:=\:\Gam(x).\]
Then for real numbers $\{a_k,c_k\}$, it is easy
to check the following limiting relation
\[\lim_{q\to1}\hyq{cccc}{q}
{q^{a_1},q^{a_2},\cdots,q^{a_n}}
{q^{c_1},q^{c_2},\cdots,q^{c_n}}_\infty
=\prod_{k=1}^n\frac{\Gam(c_k)}{\Gam(a_k)}
\quad\text{provided that}\quad
\sum_{k=1}^na_k=\sum_{k=1}^nc_k.\]
One century ago, Ramanujan~\cito{ramanujan} discovered (without proofs)
17 $\pi$-related infinite series, that were demonstrated rigorously by
Borwein brothers~\cito{bb87book} during the 1980s. Three typical formulae
are reproduced as follows:
\pnq{\frac{4}{\pi}
=&\sum_{n=0}^{\infty}
\hyp{ccc}
{\frac12,\frac12,\frac12}{\rule[1mm]{0mm}{3mm}1,\:1,\:1}_n
\frac{1+6n}{4^n}.\\
\frac{8}{\pi}
=&\sum_{n=0}^{\infty}
\hyp{ccc}
{\frac12,\frac14,\frac34}
{\rule[1mm]{0mm}{3mm}1,\:1,\:1}_n
\frac{3+20n}{(-4)^n}.\\
\frac{16}{\pi}\!
=&\sum_{n=0}^{\infty}
\hyp{ccc}
{\frac12,\frac12,\frac12}{\rule[1mm]{0mm}{3mm}1,\:1,\:1}_n
\frac{5+42n}{64^n}.}
By making use of computer algebra, Guillera~\cite{kn:gj03em,kn:gj06rj,kn:gj08rj}
found several remarkable identities including the following one
\pq{\label{g1x5hh}\frac{32}{\pi^2}=\sum_{n=0}^{\infty}
\frac{1}{16^{n}}
\hyp{c}{\frac12,\frac12,\frac12,\frac14,\frac34}{1,\,1,\,1,\,1,\,1}_n
\Big\{3+34n+120n^2\Big\}.}
Recently, there has been a growing interest in searching for $q$-analogues
of Ramanujan--like series
(cf.~\cite{kn:chu20rj,kn:chu21ijnt,kn:guo18jdea,kn:guo18itsf,kn:guo19rj,kn:gj18jdea}).
By means of the modified Abel lemma on summation by parts, the second author~\cito{chu18e}
(see also~\cito{chu20rj}) established the following $q$-analogue of \eqref{g1x5hh}:
\pq{\label{g1x5qq}\pp{c}{\hyq{c}{q}
{q^{\frac12},q^{\frac12},q^{\frac12},q^{\frac52}}{q,~q,~q,~q}_\infty
&=\sum_{n=0}^{\infty}q^{n^2}
\frac{(q^{\frac12};q)^3_n(q^{\frac32};q)_{2n}(1-q^{3n+\frac32})}
     {(q;q)^5_n(-q^{\frac12};q^{\frac12})^3_{2n+1}(1-q^{3/2})}\\
&\times\bigg\{q^n+
\frac{(1-q^n)(1-q^{3n+\frac12})(1+q^{n+\frac12})^3}
     {(1-q^{2n+\frac12})(1-q^{3n+\frac32})}
\bigg\}.}}
In 1973, Gould and Hsu~\cito{hsu} found an useful pair of inverse series
relations. Carlitz~\cito{carlitz} obtained their $q$-analogues in the
same year. These inversions together with their multiplicate forms
(cf.~\cite{kn:chu02b,kn:chu11e,kn:chu13e}) have been shown powerful in proving identities
of hypergeometric series (cf.~\cite{kn:chu93b,kn:chu94b,kn:chu94e,kn:gess82s})
and $q$-series (cf.~\cite{kn:chu95b,kn:chu09i,kn:chu02b,kn:gess83s}).
By applying the multiplicate form of the extended Carlitz inverse series
relations (cf.~\cite{kn:chu93b,kn:chu95b}) to Jackson's well--poised
$_8\phi_7$-series, we shall provide a new proof (see~Example~\ref{g1x5pp})
for identity \eqref{g1x5qq}. Numerous further $q$-series identities
will be proved likewise. These identities are significant not only
for their existence, but also for the novelty of their limiting
series as $q\to1$, namely, the classical Ramanujan--like series.
Nine remarkable ones are highlighted in advance as follows:
\pnq{\text{Example~\ref{v1x3a}:~}\quad
&\frac{\pi\Gamma^2 (\frac{1}{3}) }{6\Gamma^2(\frac{5}{6})}
=\sum_{n=0}^{\infty}
\Big(\frac{-1}{27}\Big)^n
\hyp{cccccc}{1,\:\frac{1}2,\:\frac{2}3}
{\frac43,\:\frac43,\:\frac43\rule[2mm]{0mm}{2mm}}_n
\big\{3+7 n\big\}.\\[2mm]
\text{Example~\ref{v1x3b}:~}\quad
&\frac{48\pi\Gamma^2 (\frac{2}{3}) }{\Gamma^2(\frac{1}{6})}
=\sum_{n=0}^{\infty}
\Big(\frac{-1}{27}\Big)^n
\hyp{cccccc}{1,\:\frac{1}2,\:\frac{1}3}
{\frac53,\:\frac{5}3,\:\frac{5}3\rule[2mm]{0mm}{2mm}}_n
\big\{9+28 n+21 n^2\big\}.\\[2mm]
\text{Example~\ref{v3x1a}:~}\quad
&\frac{9\sqrt3}{2^{\frac43}\pi}
=\sum_{n=0}^{\infty}
\Big(\frac{-1}{27}\Big)^n
\hyp{cccccc}{\frac13,\:\frac23,\:\frac16}
{1,\:1,\:1\rule[2mm]{0mm}{2mm}}_n
\big\{2 + 21n\big\}.\\[2mm]
\text{Example~\ref{v3x1b}:~}\quad
&\frac{27\sqrt3}{2^{\frac53}\pi}
=\sum_{n=0}^{\infty}
\Big(\frac{-1}{27}\Big)^n
\hyp{cccccc}{\frac13,\:\frac23,\:\frac56}
{1,\:1,\:1\rule[2mm]{0mm}{2mm}}_n
\big\{5 + 42n\big\}.\\[2mm]
\text{Example~\ref{v3x1c}:~}\quad
&\frac{81\sqrt3}{28\cdot2^{2/3}\pi}
=\sum_{n=0}^{\infty}
\Big(\frac{-1}{27}\Big)^n
\hyp{ccc}{\frac23,\frac73,-\frac16}
{1,~1,~2\rule[-2mm]{0mm}{2mm}}_n.\\[2mm]
\text{Example~\ref{v3x1d}:~}\quad
&\frac{3\sqrt3}{\pi\sqrt[3]4}
=\sum_{n=0}^{\infty}
\Big(\frac{-1}{27}\Big)^n
\hyp{cccccc}{\frac13,-\frac13,-\frac16}
{1,~1,~1\rule[2mm]{0mm}{2mm}}_n
\big\{1-63 n^2\big\}.\\[2mm]
\text{Example~\ref{w1+1+1a}:~}\quad
&\pp{t}{\frac{32}{\pi^2}&=\sum_{n=0}^{\infty}
\frac{(\frac12)^6_n(\frac12)^3_{2n}}{(2n)!^3(3n)!^2}
\times\frac{(1+2n)(1+4n)^2(3+8n)}{(1+3n)^2}\\
&\times\Bigg\{1+\frac{32 n (1+3 n)^2 (1+8 n)}{(1+2 n) (1+4 n)^2 (3+8 n)}
+\frac{(1+4 n) (3+4 n) (5+8 n)}{16 (2+3 n)^2 (3+8 n)}\Bigg\}.}\\[2mm]
\text{Example~\ref{w1+2d}:~}\quad
&\pp{t}{\frac{32}{\pi^2}
&=\sum_{n=0}^{\infty}(3+8n)
\dfrac{(\frac12)^3_n(\frac32)^3_{2n}(\frac32)_{3n}}
     {n!(2n)!(3n+1)!^3}\\
&\times\bigg\{1+\frac{16 n (1+3 n)^3 (1+10 n)}{(1+4 n)^3 (1+6 n) (3+8 n)}
+\frac{3 (1+2 n)^3 (7+10 n)}{16 (2+3 n)^3 (3+8 n)}\Big\}.}\\[2mm]
\text{Example~\ref{w1+2e}:~}\quad
&\pp{t}{\frac{9\pi^2}{16}
&=\sum_{n=0}^{\infty}
\dfrac{(2n)!(4n+1)!}
     {(6n+1)!(1+2n)(1+3n)}\\
&\times\bigg\{1+\frac{9 (1+2 n) (1+3 n) (1+5 n)}{2 (1+4 n)^2 (1+6 n)}
+\frac{3(1+n) (1+2 n)^2 (4+5 n)}{2 (2+3 n) (3+4 n) (5+6 n)^2)}\bigg\}.}}
The rest of the paper will be organized as follows.
By applying the multiplicate form of the extended Carlitz inverse
series relations to Jackson's $q$-analogue of Dougall's well--poised
$_7F_6$-sum, we shall derive, in the next section, two general dual
formulae (Theorems~\ref{thm=U} and ~\ref{thm=V}). Then under the
partition pattern $n=\pav{\frac{n}2}+\pav{\frac{n+1}2}$, we shall deduce,
in Section~3, from the duplicate inversions a number of $q$-series
identities whose limiting cases as $q\to1$ result in classical
$\pi$-related Ramanujan--like series of convergence rate
``$\frac1{16}$". The triplicate inversions under the partition pattern
$n=\pav{\frac{n}3}+\pav{\frac{n+1}3}+\pav{\frac{n+2}3}$ will be utilized
in Section~4 to establish further formulae for $q$-series, where the
``reverse bisection method" is employed to reduce the resulting series
dramatically to simpler but more elegant ones.

\section{Inversion Approach to Jackson's Well--Poised $_8\phi_7$-Series}
In 1973, Gould and Hsu~\cito{hsu} discovered a useful pair
of inverse series relations, which can equivalently be
reproduced below. Let $\{a_i,b_i\}$ be  two complex
sequences such that the $\vph$-polynomials defined by
\pq{\vph(x;0)\equiv1\label{phi-pp}
\quad\text{and}\quad
\vph(x;n)=\prod_{k=0}^{n-1}(a_k+xb_k)
\quad\text{for}\quad n\in\mb{N}}
differ from zero for $x,\:n\in\mb{N}_0$.
Then there hold the inverse series relations
\pmq{f(n)&=\label{inv+gh}
\sum_{k=0}^{n}
(-1)^k\binm{n}{k}
\vph(k;n)\:g(k),\\
g(n)&=\label{inv-gh}
\sum_{k=0}^{n}
(-1)^k\binm{n}{k}
\frac{a_k+kb_k}{\vph(n;k+1)}
\:f(k).}
Define the Gaussian binomial coefficient by
\[\binq{m}{n}=\frac{(q;q)_m}{(q;q)_n(q;q)_{m-n}}
=\frac{(q^{m-n+1};q)_n}{(q;q)_n},
\quad\text{where}\quad m,\:n\in\mb{N}.\]
Carlitz~\cito{carlitz} found the following $q$-analogues
\pmq{f(n)&=\label{carlitz+f}
\sum_{k=0}^n(-1)^k
\binq{n}{k}\vph(q^{-k};n)g(k),\\
g(n)&=\label{carlitz+g}
\sum_{k=0}^n(-1)^k\binq{n}{k}q^{\binm{n-k}{2}}
\frac{a_k+q^{-k}b_k}{\vph(q^{-n};k+1)}f(k);}
which were utilized to derive $q$-series identities
\cite{kn:chu94d,kn:chu95b,kn:gess83s}.
There is an extended form (cf.~\cite{kn:chu93b,kn:chu95b}):
\pmq{f(n)&=\label{inv+v}
\sum_{k=0}^{n}(-1)^k\binq{n}{k}
\vph(q^k\sig;n)\vph(q^{-k};n)
\frac{1-q^{2k}\sig}{(q^n\sig;q)_{k+1}}g(k),\\
g(n)&=\label{inv-v}
\sum_{k=0}^{n}
(-1)^k\binq{n}{k}q^{\binm{n-k}{2}}
\frac{(a_k+q^k\sig b_k)(a_k+q^{-k}b_k)}
{\vph(q^n\sig;k+1)\vph(q^{-n};k+1)}(q^k\sig;q)_n
\:f(k).}

In order to facilitate the subsequent application,
we can reformulate the above inversions as
\pmq{
F(n)&=\label{inv+w}
\sum_{k=0}^{n}
(-1)^k\binq{n}{k}
\frac{\phi(q^k;n)}{\psi(q^n;k+1)}
\frac{\psi(q^k;k+1)}{\psi(q^k;k)}G(k),\\
G(n)&=\label{inv-w}
\sum_{k=0}^{n}
(-1)^k\binq{n}{k}q^{\binm{n-k}{2}}
\frac{\psi(q^k;n)}{\phi(q^n;k+1)}
\frac{\phi(q^k;k+1)}{\phi(q^k;k)}\:F(k);}
where the polynomial sequences $\{\phi(x;n),\psi(x;n)\}$
are defined respectively by
\pq{\phi(x;n)=\vph(x\sig;n)\vph(x^{-1};n)
\quad\text{and}\quad
\psi(x;n)=(x\sig;q)_n.}

For the fundamental identity of Dougall~\cito{dougall}
about the very well--poised terminating $_7F_6$-series,
Jackson~\cito{jackson} (cf.~Bailey~\citu{bailey}{\S8.3})
discovered, one century ago, the following $q$-analogue:
\pq{\label{jackson}
\pp{c}{\Ome_n(a;b,c,d):=&\sum_{k=0}^n\frac{1-q^{2k}a}{1-a}
\hyq{ccc}{q}{a,~b,~c,~d,~e,~q^{-n}}
{q,qa/b,qa/c,qa/d,qa/e,q^{n+1}a}_kq^k\\
=&\hyq{c}{q}{qa,qa/bc,qa/bd,qa/cd}{qa/b,qa/c,qa/d,qa/bcd}_n,
\quad\text{where}\quad q^{n+1}a^2=bcde.}}


For a real number $x$, denote by $\pav{x}$ 
the greatest integer $\le x$. 
Suppose that for all $n\in\mb{N}_0$, there exist eight integers
$\{\veps_b,\veps_c,\veps_d,\veps_e,\lam_b,\lam_c,\lam_d,\lam_e\}\subset\mb{N}_0$
with $\boxed{\Lam:=\lam_b+\lam_c+\lam_d+\lam_e}$ such that
\pq{\label{p-pattern}
n=\pavv{\frac{\veps_b+n\lam_b}{\Lam}}
+\pavv{\frac{\veps_c+n\lam_c}{\Lam}}
+\pavv{\frac{\veps_d+n\lam_d}{\Lam}}
+\pavv{\frac{\veps_e+n\lam_e}{\Lam}}.}
For the sake of brevity, denote this partition by
\[n=\ang{n}_b+\ang{n}_c+\ang{n}_d+\ang{n}_e,\]
where
\pnq{\ang{n}_b&:=\pavv{\frac{\veps_b+n\lam_b}{\Lam}},\\
\ang{n}_c&:=\pavv{\frac{\veps_c+n\lam_c}{\Lam}},\\
\ang{n}_d&:=\pavv{\frac{\veps_d+n\lam_d}{\Lam}},\\
\ang{n}_e&:=\pavv{\frac{\veps_e+n\lam_e}{\Lam}}.}

Then Jackson's formula \eqref{jackson} can be reformulated as
\pnq{\Ome_n\Big(a;bq^{\ang{n}_b},cq^{\ang{n}_c},dq^{\ang{n}_d})
&=\frac{(qa;q)_n(bc/a;q)_{\ang{n}_b+\ang{n}_c}(qa/bc;q)_{\ang{n}_d+\ang{n}_e}}
{(b/a;q)_{\ang{n}_b}(c/a;q)_{\ang{n}_c}(d/a;q)_{\ang{n}_d}(qa/bcd;q)_{\ang{n}_e}}\\
&\times\frac{(bd/a;q)_{\ang{n}_b+\ang{n}_d}(qa/bd;q)_{\ang{n}_c+\ang{n}_e}}
{(qa/b;q)_{\ang{n}_c+\ang{n}_d+\ang{n}_e}(qa/c;q)_{\ang{n}_b+\ang{n}_d+\ang{n}_e}}\\
&\times\frac{(cd/a;q)_{\ang{n}_c+\ang{n}_d}(qa/cd;q)_{\ang{n}_b+\ang{n}_e}}
{(qa/d;q)_{\ang{n}_b+\ang{n}_c+\ang{n}_e}(bcd/a;q)_{\ang{n}_b+\ang{n}_c+\ang{n}_d}}.}

In view of the relation (which is also true when ``$b$" is replaced by one of ``$c,d,e$")
\[\frac{(bq^{\ang{n}_b};q)_k}{(q^{1-\ang{n}_b}a/b;q)_k}
=q^{k\ang{n}_b}\frac{(b;q)_k(q^kb;q)_{\ang{n}_b}(q^{-k}b/a;q)_{\ang{n}_b}}
{(qa/b;q)_k(b;q)_{\ang{n}_b}(b/a;q)_{\ang{n}_b}},\]
it is not difficult to check that Jackson's formula is equivalent
to the following $q$-binomial sum
\pnq{\sum_{k=0}^n&(-1)^k\binq{n}{k}\frac{1-q^{2k}a}{(q^na;q)_{k+1}}
\hyq{c}{q}{a,~b,~c,~d,~qa^2/bcd}{qa/b,qa/c,qa/d,bcd/a}_k
q^{\binm{k+1}{2}}\\
&~\times
(q^kb;q)_{\ang{n}_b}(q^kc;q)_{\ang{n}_c}(q^kd;q)_{\ang{n}_d}(q^{1+k}a^2/bcd;q)_{\ang{n}_e}\\
&~\times
(q^{-k}b/a;q)_{\ang{n}_b}(q^{-k}c/a;q)_{\ang{n}_c}(q^{-k}d/a;q)_{\ang{n}_d}(q^{1-k}a/bcd;q)_{\ang{n}_e}\\
=~&{(a;q)_n(b;q)_{\ang{n}_b}(c;q)_{\ang{n}_c}(d;q)_{\ang{n}_d}(qa^2/bcd;q)_{\ang{n}_e}}\\[1mm]
\times~&
\frac{(bc/a;q)_{\ang{n}_b+\ang{n}_c}(bd/a;q)_{\ang{n}_b+\ang{n}_d}(cd/a;q)_{\ang{n}_c+\ang{n}_d}}
{(qa/b;q)_{\ang{n}_c+\ang{n}_d+\ang{n}_e}(qa/c;q)_{\ang{n}_b+\ang{n}_d+\ang{n}_e}}\\
\times~&
\frac{(qa/bc;q)_{\ang{n}_d+\ang{n}_e}(qa/bd;q)_{\ang{n}_c+\ang{n}_e}(qa/cd;q)_{\ang{n}_b+\ang{n}_e}}
     {(qa/d;q)_{\ang{n}_b+\ang{n}_c+\ang{n}_e}(bcd/a;q)_{\ang{n}_b+\ang{n}_c+\ang{n}_d}}.}

Observing that this equation matches exactly
to \eqref{inv+w} specified by ``$\sig=a$" and
\pnq{
\phi(x;n)&=(bx;q)_{\ang{n}_b}(cx;q)_{\ang{n}_c}(dx;q)_{\ang{n}_d}(qa^2x/bcd;q)_{\ang{n}_e}\\
&\times(b/ax;q)_{\ang{n}_b}(c/ax;q)_{\ang{n}_c}(d/ax;q)_{\ang{n}_d}(qa/bcdx;q)_{\ang{n}_e},\\
\psi(x;n)&=(ax;q)_n;}
as well as
\pnq{F(n)&={(a;q)_n(b;q)_{\ang{n}_b}(c;q)_{\ang{n}_c}(d;q)_{\ang{n}_d}(qa^2/bcd;q)_{\ang{n}_e}}\\[1mm]
&\times~
\frac{(bc/a;q)_{\ang{n}_b+\ang{n}_c}(bd/a;q)_{\ang{n}_b+\ang{n}_d}(cd/a;q)_{\ang{n}_c+\ang{n}_d}}
{(qa/b;q)_{\ang{n}_c+\ang{n}_d+\ang{n}_e}(qa/c;q)_{\ang{n}_b+\ang{n}_d+\ang{n}_e}}\\
&\times~
\frac{(qa/bc;q)_{\ang{n}_d+\ang{n}_e}(qa/bd;q)_{\ang{n}_c+\ang{n}_e}(qa/cd;q)_{\ang{n}_b+\ang{n}_e}}
     {(qa/d;q)_{\ang{n}_b+\ang{n}_c+\ang{n}_e}(bcd/a;q)_{\ang{n}_b+\ang{n}_c+\ang{n}_d}},\\
G(n)&=\hyq{c}{q}{a,~b,~c,~d,~qa^2/bcd}{qa/b,qa/c,qa/d,bcd/a}_n
q^{\binm{n+1}{2}};}
we may state the dual relation corresponding to \eqref{inv-w} as
\[G(n)=\sum_{k=0}^{n}
(-1)^k\binq{n}{k}q^{\binm{n-k}{2}}
\frac{\psi(q^k;n)}{\phi(q^n;k+1)}
\frac{\phi(q^k;k+1)}{\phi(q^k;k)}F(k).\]
This can explicitly be highlighted, as the following lemma.
\begin{lemm}\label{mm=H}
The following identity holds
\pnq{\hyq{c}{q}{b,~c,~d,~qa^2/bcd}{qa/b,qa/c,qa/d,bcd/a}_n
=&\sum_{k=0}^nq^{k-n}
\frac{(q^na;q)_k(q^{-n};q)_k}{\phi(q^n;k+1)}
\frac{\phi(q^k;k+1)}{\phi(q^k;k)}H(k),}
where $H(n)$ is the quotient of shifted factorials given by
\pnq{H(n)&=\frac{(b;q)_{\ang{n}_b}(c;q)_{\ang{n}_c}(d;q)_{\ang{n}_d}(qa^2/bcd;q)_{\ang{n}_e}}
{(q;q)_n~(bcd/a;q)_{\ang{n}_b+\ang{n}_c+\ang{n}_d}}\\[1mm]
&\times~
\frac{(bc/a;q)_{\ang{n}_b+\ang{n}_c}(qa/bc;q)_{\ang{n}_d+\ang{n}_e}}
{(qa/b;q)_{\ang{n}_c+\ang{n}_d+\ang{n}_e}}\\
&\times~
\frac{(bd/a;q)_{\ang{n}_b+\ang{n}_d}(qa/bd;q)_{\ang{n}_c+\ang{n}_e}}
     {(qa/c;q)_{\ang{n}_b+\ang{n}_d+\ang{n}_e}}\\
&\times~
\frac{(cd/a;q)_{\ang{n}_c+\ang{n}_d}(qa/cd;q)_{\ang{n}_b+\ang{n}_e}}
{(qa/d;q)_{\ang{n}_b+\ang{n}_c+\ang{n}_e}}.}
\end{lemm}

Recall that for all $k\in\mb{N}_0$, we have
\[\ang{k+1}_b+\ang{k+1}_c+\ang{k+1}_d+\ang{k+1}_e=k+1\]
and
\pnq{\phi(q^n;k+1)=(q^nb;q)_{\ang{k+1}_b}(q^nc;q)_{\ang{k+1}_c}(q^nd;q)_{\ang{k+1}_d}(q^{1+n}a^2/bcd;q)_{\ang{k+1}_e}&\\
\times(q^{-n}b/a;q)_{\ang{k+1}_b}(q^{-n}c/a;q)_{\ang{k+1}_c}(q^{-n}d/a;q)_{\ang{k+1}_d}(q^{1-n}a/bcd;q)_{\ang{k+1}_e}.&}
When $n\to\infty$, it is not hard to check the limiting relation
\pnq{q^{k-n}
\frac{(q^na;q)_k(q^{-n};q)_k}{\phi(q^n;k+1)}
&\Rightarrow-q^{\binm{k+1}2-\binm{\ang{k+1}_b}2-\binm{\ang{k+1}_c}2-\binm{\ang{k+1}_d}2-\binm{\ang{k+1}_e}2}\\
&\times\Big(\frac{a}{b}\Big)^{\ang{k+1}_b}\Big(\frac{a}{c}\Big)^{\ang{k+1}_c}
\Big(\frac{a}{d}\Big)^{\ang{k+1}_d}\Big(\frac{bcd}{qa}\Big)^{\ang{k+1}_e}.}

Now let $n\to\infty$ in Lemma~\ref{mm=H}
and then determine the limit of the right member through the Weierstrass
$M$-test on uniformly convergent series (cf. Stromberg~\citu{karl}{\S3.106}).
After some routine simplification, the resulting limiting relation may be
stated as in the following proposition.
\begin{prop} The following identity holds\label{pp=H}
\[\hyq{c}{q}{b,~c,~d,~qa^2/bcd}{qa/b,qa/c,qa/d,bcd/a}_{\infty}
=\sum_{k=0}^{\infty}\Theta(k+1)
\bigg\{\frac{\phi(q^k;k+1)}{-\phi(q^k;k)}\bigg\}H(k),\]
where the monomial $\Theta(m)$ is defined by
\pnq{\Theta(m)&=
q^{\binm{m}2-\binm{\ang{m}_b}2-\binm{\ang{m}_c}2
	-\binm{\ang{m}_d}2-\binm{\ang{m}_e}2}\\
&\times~
\Big(\frac{a}{b}\Big)^{\ang{m}_b}\Big(\frac{a}{c}\Big)^{\ang{m}_c}
\Big(\frac{a}{d}\Big)^{\ang{m}_d}\Big(\frac{bcd}{qa}\Big)^{\ang{m}_e}.}
\end{prop}

To facilitate understanding, we explicitly display the expression from this proposition:
\pnq{&\xpuad\hyq{c}{q}{b,~c,~d,~qa^2/bcd}{qa/b,qa/c,qa/d,bcd/a}_\infty
=\sum_{k=0}^{\infty}
\bigg\{\frac{\phi(q^k;k+1)}{-\phi(q^k;k)}\bigg\}\\
&\times~
q^{\binm{k+1}2-\binm{\ang{k+1}_b}2-\binm{\ang{k+1}_c}2-\binm{\ang{k+1}_d}2-\binm{\ang{k+1}_e}2}\\
&\times\Big(\frac{a}{b}\Big)^{\ang{k+1}_b}\Big(\frac{a}{c}\Big)^{\ang{k+1}_c}
\Big(\frac{a}{d}\Big)^{\ang{k+1}_d}\Big(\frac{bcd}{qa}\Big)^{\ang{k+1}_e}\\
&\times~
\frac{(b;q)_{\ang{k}_b}(c;q)_{\ang{k}_c}(d;q)_{\ang{k}_d}(qa^2/bcd;q)_{\ang{k}_e}}
     {(q;q)_k~(bcd/a;q)_{\ang{k}_b+\ang{k}_c+\ang{k}_d}}\\[1mm]
&\times~
\frac{(bc/a;q)_{\ang{k}_b+\ang{k}_c}(qa/bc;q)_{\ang{k}_d+\ang{k}_e}}
     {(qa/b;q)_{\ang{k}_c+\ang{k}_d+\ang{k}_e}}\\
&\times~
\frac{(bd/a;q)_{\ang{k}_b+\ang{k}_d}(qa/bd;q)_{\ang{k}_c+\ang{k}_e}}
     {(qa/c;q)_{\ang{k}_b+\ang{k}_d+\ang{k}_e}}\\
&\times~
\frac{(cd/a;q)_{\ang{k}_c+\ang{k}_d}(qa/cd;q)_{\ang{k}_b+\ang{k}_e}}
     {(qa/d;q)_{\ang{k}_b+\ang{k}_c+\ang{k}_e}}.}


Classifying the summation indices $\boxed{k\to k+n\Lam}$ with
respect to their residues modulo $\Lam$, we can split the sum
in Proposition~\ref{pp=H} into $\Lam$ partial sums:
\pnq{&\hyq{c}{q}{b,~c,~d,~qa^2/bcd}{qa/b,qa/c,qa/d,bcd/a}_\infty\\
=&\sum_{n=0}^{\infty}\sum_{k=0}^{\Lam-1}\Theta(1+k+n\Lam)
\bigg\{\frac{\phi(q^{k+n\Lam};{1+k+n\Lam})}
      {-\phi(q^{k+n\Lam};{k+n\Lam})}\bigg\}
H(k+n\Lam)\\
=&\sum_{n=0}^{\infty}H(n\Lam)\Theta(n\Lam)
\sum_{k=0}^{\Lam-1}\frac{\Theta(1+k+n\Lam)}{\Theta(n\Lam)}
\bigg\{\frac{\phi(q^{k+n\Lam};{1+k+n\Lam})}{-\phi(q^{k+n\Lam};{k+n\Lam})}\bigg\}
\frac{H(k+n\Lam)}{H(n\Lam)}.}

In general, by singling out the initial $\del$-terms, where $\boxed{0\le\del<\Lam}$,
we have two further expressions, that will efficiently be utilized
to derive concrete $q$-series identities in the next two sections.
\begin{thm}[Nonterminating series identity]\label{thm=U}
\pnq{\hyq{c}{q}{b,~c,~d,~qa^2/bcd}{qa/b,qa/c,qa/d,bcd/a}_{\infty}
=\sum_{k=0}^{\del-1}\Theta(k+1)
\bigg\{\frac{\phi(q^k;k+1)}{-\phi(q^k;k)}\bigg\}&H(k)\\
+\sum_{n=0}^{\infty}H(n\Lam)\Theta(n\Lam)
\sum_{i=\del}^{\del-1+\Lam}\frac{\Theta(1+i+n\Lam)}{\Theta(n\Lam)}
\bigg\{\frac{\phi(q^{i+n\Lam};1+i+n\Lam)}{-\phi(q^{i+n\Lam};i+n\Lam)}\bigg\}
&\frac{H(i+n\Lam)}{H(n\Lam)}.}
\end{thm}
\begin{thm}[Nonterminating series identity]\label{thm=V}
\pnq{\hyq{c}{q}{b,~c,~d,~qa^2/bcd}{qa/b,qa/c,qa/d,bcd/a}_{\infty}
=\sum_{k=0}^{\del-1}\Theta(k+1)
\bigg\{\frac{\phi(q^k;k+1)}{-\phi(q^k;k)}\bigg\}&H(k)\\
+\sum_{n=1}^{\infty}H(n\Lam)\Theta(n\Lam)
\!\sum_{j=\del-\Lam}^{\del-1}\puab\frac{\Theta(1+j+n\Lam)}{\Theta(n\Lam)}
\bigg\{\frac{\phi(q^{j+n\Lam};1+j+n\Lam)}{-\phi(q^{j+n\Lam};j+n\Lam)}\bigg\}
&\frac{H(j+n\Lam)}{H(n\Lam)}.}
\end{thm}

\section{Series from Duplicate Inversions}
For all $n\in\mb{N}_0$, it is well known that
$\boxed{n=\pavv{\tfrac{n}2}+\pavv{\tfrac{1+n}2}}$.
Then it is not difficult to check that Jackson's formula \eqref{jackson}
is equivalent to the following one
\[\Ome_n\big(a;q^{\pav{\tfrac{n}2}}b,c,q^{\pav{\tfrac{1+n}2}}d\big)
=\hyq{c}{q}{qa,\quad bd/a}{qa/c,bcd/a}_n
\hyq{c}{q}{qa/cd,bc/a}{qa/d,b/a}_{\pavv{\frac{n}2}}
\hyq{c}{q}{qa/bc,cd/a}{qa/b,d/a}_{\pavv{\frac{1+n}2}}\]
with its parameters subject to $\boxed{qa^2=bcde}$.
In this case, the partition pattern defined by \eqref{p-pattern}
can be explicitly determined as
\pq{\label{p2pattern}\boxed{\Lam=2:\mult{cccc}
{\veps_b=0,&\veps_c=0,&\veps_d=1,&\veps_e=0;\\
\lam_b=1,&\lam_c=0,&\lam_d=1,&\lam_e=0.}}}

\subsection{} \
Letting $\boxed{\del=0}$ in Theorem~\ref{thm=U}, we derive, under the
above partition pattern, the following quite general formula involving
four free parameters $\{a,b,c,d\}$.
\begin{thm}[Nonterminating series identity]\label{thm=2U}
\pnq{\hyq{c}{q}{b,~c,~qd,~qa^2/bcd}{qa/b,qa/c,a/d,bcd/a}_\infty
=\sum_{n=0}^{\infty}W_n(a,b,c,d)
\Big(\frac{q^{n}a^2}{bd}\Big)^n&\\
\times\frac{\ff{bd/a;q}{2n}\ff{b,qd,q^2a/bc,bc/a,qa/cd,cd/a;q}{n}}
{\ff{q,qa/c,bcd/a;q}{2n}\ff{qa/b,qa/d;q}{n}},&}
where the weight function is given by
\pnq{W_n(a,b,c,d)
&=\frac{(1-qa/bc)(1-q^ncd/a)(1-q^{2n}bd/a)(1-q^{1+3n}b)}
{(1-a/d)(1-q^{1+2n})(1-q^{1+2n}a/c)(1-q^{2n}bcd/a)}\\
\times\bigg\{q^n\tfrac{a}{d}
&+\tfrac{(1-q^na/d)(1-q^{1+2n})(1-q^{1+2n}a/c)(1-q^{2n}bcd/a)(1-q^{3n}d)}
{(1-q^{1+n}a/bc)(1-q^nd)(1-q^ncd/a)(1-q^{2n}bd/a)(1-q^{1+3n}b)}\bigg\}.}
\end{thm}

Letting $\boxed{a=b=c=d=q^{\frac12}}$ in this theorem and then simplifying
slightly the resulting expression, we get immediately the following infinite
series identity, that was derived by the second author~\cito{chu18e}
by making use of the modified Abel lemma on summation by parts.
\begin{exam}[$\boxed{a=b=c=d=q^{\frac12}}$ in Theorem~\ref{thm=2U}]\label{g1x5pp}
\pnq{\hyq{c}{q}{q^{\frac12},q^{\frac12},q^{\frac12},q^{\frac52}}{q,~q,~q,~q}_\infty
&=\sum_{n=0}^{\infty}q^{n^2}
\frac{(q^{\frac12};q)^3_n(q^{\frac32};q)_{2n}(1-q^{3n+\frac32})}
     {(q;q)^5_n(-q^{\frac12};q^{\frac12})^3_{2n+1}(1-q^{3/2})}\\
&\times\bigg\{q^n+
\frac{(1-q^n)(1-q^{3n+\frac12})(1+q^{n+\frac12})^3}
     {(1-q^{2n+\frac12})(1-q^{3n+\frac32})}
\bigg\}.}
\end{exam}
When $q\to1$, we recover the following beautiful series discovered 
by Guillera~\cite{kn:gj03em,kn:gj06rj}: 
\[\frac{32}{\pi^2}
=\sum_{n=0}^{\infty}
\frac{1}{16^{n}}
\hyp{c}{\frac12,\frac12,\frac12,\frac14,\frac34}{1,\,1,\,1,\,1,\,1}_n
\Big\{3+34n+120n^2\Big\}.\]

Further $q$-analogues of Ramanujan--like $\pi$-related series can be derived
similarly, where the parameter settings for $\{a,b,c,d\}$ are highlighted
in the headers for all the examples.
\begin{exam}[$\boxed{a=q^{\frac32},~b=q^{\frac12},~c=q^{\frac32},~d=q^{\frac12}}$ in Theorem~\ref{thm=2U}]
\pnq{\hyq{c}{q}{q^{\frac32},q^{\frac32},q^{\frac32},q^{\frac32}}{q,~q,~q^2,~q^2}_\infty
=\sum_{n=0}^{\infty}q^{n^2+2n}
\frac{(q^{\frac12};q)^3_n(q^{-\frac12};q)_{2n}(1+q^{1/2})}
     {(q^2;q^2)^3_{n}(q;q)_{n}(q^2;q)_{n}(-q^{\frac12};q)^3_{n}}&\\
\times\frac{1-q^{3n+\frac12}}{1-q^{1/2}}
\bigg\{1+
\frac{q^{1+n}(1-q^{2n-\frac12})(1-q^{3n+\frac32})}
{(1-q^{1+n})(1-q^{3n+\frac12})(1+q^{n+\frac12})^3}
\bigg\}&.}
\end{exam}
When $q\to1$, we recover the following known series 
(cf.~\cite{kn:chu14mc}): 
\[\frac{128}{\pi^2}
=\sum_{n=0}^{\infty}
\frac{1}{16^{n}}
\hyp{c}{\frac12,\frac12,\frac12,\frac14,-\frac14}{1,\,1,\,1,\,2,\,2}_n
\Big\{13+118n+120n^2\Big\}.\]

\begin{exam}[$\boxed{a=q^{\frac32},~b=q^{\frac12},~c=q^{\frac12},~d=q^{\frac32}}$ in Theorem~\ref{thm=2U}]
\pnq{\hyq{c}{q}
{q^{\frac32},q^{\frac52},q^{\frac52},q^{\frac52}}{q^2,~q^2,~q^2,~q^3}_\infty
&=\sum_{n=0}^{\infty}q^{n^2+n}
\frac{(q^{\frac12};q)^2_n(q^{-\frac12};q)_n(q^{\frac32};q)_n(q^{\frac52};q)^2_n(q^{\frac32};q)_{2n}}
     {(q^2;q)^2_{2n}(q^3;q)_{2n}(q;q)_{n}(q^2;q)_{n}}\\
&\times\frac{1-q^{3n+\frac32}}{1-q^{3/2}}
\bigg\{q^{n}+\frac{(1-q^{n})(1-q^{2n+1})^2(1-q^{2n+2})}
     {(1-q^{n+\frac12})(1-q^{n+\frac32})^2(1-q^{2n+\frac12})}
\bigg\}.}
\end{exam}
This provides a $q$-analogue for the following 
classical series:  
\[ \frac{256}{3 \pi ^2}
=\sum_{n=0}^{\infty}
\hyp{c}
{\frac{1}{2},-\frac{1}{2},\frac{3}{2},\frac{1}{4},\frac{3}{4}}
{\\[-3mm] 1,1,1,~2,~2}_n
\frac{80 n^3+148 n^2+80 n+9}{16^n}.\]

\begin{exam}[$\boxed{a=q^{\frac32},~b=q^{\frac32},~c=q^{\frac12},~d=q^{\frac32}}$ in Theorem~\ref{thm=2U}]
\pnq{\hyq{c}{q}{q^{\frac32},q^{\frac32},q^{\frac52},q^{\frac72}}{q,~q^2,~q^3,~q^3}_\infty
&=\sum_{n=0}^{\infty}q^{n^2}
\frac{(q^{\frac12};q)^2_n(q^{\frac52};q)_n(q^{\frac52};q)_{2n}(1-q^{3n+\frac52})}
     {(q;q)^2_{n}(q^2;q^2)_{n}(q^4;q^2)^2_{n}(-q^{\frac32};q)^3_{n}(1-q^{5/2})}\\
&\times\bigg\{q^{n}+
\frac{(1-q^{n})(1-q^{2n+1})(1-q^{2n+2})^2(1-q^{3n+\frac32})}
     {(1-q^{n+\frac12})^2(1-q^{n+\frac32})(1-q^{2n+\frac32})(1-q^{3n+\frac52})}
\bigg\}.}
\end{exam}
This provides a $q$-analogue for the following 
classical series:   
\[ \frac{512}{\pi ^2}
=\sum_{n=0}^{\infty}
\hyp{c}
{ \frac{1}{2},\frac{1}{2},\frac{3}{2},\frac{3}{4},\frac{5}{4}}
{\\[-3mm] 1,1,1,2,2}_n
\frac{240 n^3+532 n^2+336 n+45}{16^n}.\]

\begin{exam}[$\boxed{a=q^{-\frac12},~b=q^{-\frac56},~c=q^{\frac16},~d=q^{\frac16}}$ in Theorem~\ref{thm=2U}]
\pnq{\hyq{c}{q}{q^{\frac12},q^{\frac76},q^{\frac76},q^{\frac{13}6}}
               {q^{2},q^{\frac13},q^{\frac43},q^{\frac{4}3}}_\infty
=\sum_{n=0}^{\infty}q^{n^2+\frac23n}
\frac{\ff{q^{\frac16},q^{-\frac16},q^{-\frac{5}6},q^{\frac{7}6},
         q^{\frac{11}6},q^{\frac{13}6};q}{n}\ff{q^{\frac56};q}{2n}}
     {(q^{\frac{1}3};q)_{n}(q^{\frac{4}3};q)_{n}(q;q)_{2n}(q^{2};q)_{2n}(q^{\frac{4}3};q)_{2n}}&\\
\times\frac{1-q^{3n+\frac16}}
     {1-q^{\frac16}}
\bigg\{1-
\frac{(1-q^{\frac23-n})(1-q^{2n+\frac13})(1-q^{2n})(1-q^{2n+1})}
     {(1-q^{n+\frac16})(1-q^{n+\frac56})(1-q^{n+\frac{7}6})(1-q^{2n-\frac16})}
\bigg\}&.}
\end{exam}
This provides a $q$-analogue for the following 
classical series:  
\[\frac{180 \Gamma (\frac{2}{3})^3}{\pi ^2}
=\sum_{n=0}^{\infty}
\hyp{c}
{ \frac{1}{6},-\frac{1}{6},\frac{5}{6},-\frac{5}{6},
\frac{7}{6},-\frac{1}{12},\frac{5}{12},\frac{19}{18}}
{\\[-3mm] 1,~ 1,~ \frac{1}{2},~\frac{3}{2},~\frac{1}{3},~
\frac{2}{3},~\frac{4}{3},~\frac{1}{18}}_n
\frac{35+228 n-540 n^2-2160 n^3}{16^n}.\]

\begin{exam}[$\boxed{a=q^{\frac12},~b=q^{\frac16},~c=q^{-\frac56},~d=q^{\frac76}}$ in Theorem~\ref{thm=2U}]
\pnq{\hyq{c}{q}{q^{\frac32},q^{\frac16},q^{\frac{13}6},q^{\frac{19}6}}
               {q^{2},q^{\frac13},q^{\frac43},q^{\frac{10}3}}_\infty
=\sum_{n=0}^{\infty}q^{n^2+\frac23n}
\frac{\ff{q^{\frac16},q^{\frac56},q^{\frac{7}6},q^{-\frac{7}6},
         q^{\frac{13}6},q^{\frac{19}6};q}{n}(q^{\frac{11}6};q)_{2n}}
     {(q^{\frac{1}3};q)_{n}(q^{\frac{4}3};q)_{n}(q;q)_{2n}(q^{2};q)_{2n}(q^{\frac{10}3};q)_{2n}}&\\
\times\frac{1-q^{3n+\frac76}}
     {1-q^{\frac76}}
\bigg\{1-
\frac{(1-q^{\frac23-n})(1-q^{2n+\frac73})(1-q^{2n})(1-q^{2n+1})}
     {(1-q^{n-\frac16})(1-q^{n+\frac76})(1-q^{n+\frac{13}6})(1-q^{2n+\frac56})}
\bigg\}&.}
\end{exam}
This provides a $q$-analogue for the following 
classical series:  
\[ \frac{960 \Gamma (\frac{2}{3})^3}{7 \pi ^2}
=\sum_{n=0}^{\infty}
\hyp{c}
{ \frac{1}{6},-\frac{1}{6},\frac{7}{6},-\frac{7}{6},
\frac{13}{6},\frac{5}{12},\frac{11}{12},\frac{25}{18}}
{\\[-3mm] 1,~ 1,~ \frac{1}{2},~\frac{3}{2},~\frac{1}{3},~
\frac{4}{3},~\frac{5}{3},~\frac{7}{18}}_n
\frac{65+372 n-756 n^2-2160 n^3}{16^n}.\]

\begin{exam}[$\boxed{a=q^{\frac12},~b=q^{\frac16},~c=q^{\frac76},~d=q^{\frac76}}$ in Theorem~\ref{thm=2U}]
\pnq{\hyq{c}{q}{q^{-\frac12},q^{\frac16},q^{\frac76},q^{\frac{13}6}}
               {q^{2},q^{\frac13},q^{-\frac23},q^{\frac{4}3}}_\infty
=\sum_{n=0}^{\infty}q^{n^2-\frac{n}3}
\frac{\ff{q^{\frac16},q^{\frac16},q^{\frac{5}6},q^{-\frac56},
         q^{\frac76},q^{\frac{11}6};q}{n}(q^{\frac56};q)_{2n}}
     {(q^{-\frac23};q)_{n}(q^{\frac43};q)_{n}(q;q)_{2n}(q^{2};q)_{2n}(q^{\frac{1}3};q)_{2n}}&\\
\times\frac{1-q^{3n+\frac76}}
     {1-q^{\frac76}}
\bigg\{1+q^{n-\frac23}
\frac{(1-q^{n+\frac16})(1-q^{n+\frac76})(1-q^{n+\frac{11}6})(1-q^{2n+\frac56})}
     {(1-q^{n-\frac23})(1-q^{2n+\frac13})(1-q^{2n+1})(1-q^{2n+2})}
\bigg\}&.}
\end{exam}
This provides a $q$-analogue for the following 
classical series:  
\[ \frac{7776 \Gamma (\frac{2}{3})^3}{7 \pi ^2}
=\sum_{n=0}^{\infty}
\hyp{c}
{\frac{1}{6},\frac{5}{6},\frac{5}{6},-\frac{5}{6},\frac{7}{6},
   \frac{11}{6},\frac{5}{12},\frac{11}{12},\frac{25}{18}}
{\\[-3mm] 1,2,~\frac{3}{2},\frac{3}{2},~\frac{1}{3},
\frac{2}{3},\frac{4}{3},-\frac{1}{6},~\frac{7}{18}}_n
\frac{360 n^2+546 n+191}{16^n}.\]

\begin{exam}[$\boxed{a=q^{\frac12},~b=q^{\frac56},~c=q^{-\frac16},~d=q^{-\frac16}}$ in Theorem~\ref{thm=2U}]
\pnq{\hyq{c}{q}{q^{\frac32},q^{\frac56},q^{\frac{11}6},q^{\frac{17}6}}
               {q^{2},q^{\frac23},\:q^{\frac53},\:q^{\frac{8}3}}_\infty
=\sum_{n=0}^{\infty}q^{n^2+\frac{4}3n}
\frac{(q^{\frac16};q)^2_n(q^{\frac56};q)^2_n(q^{\frac{11}6};q)^2_n(q^{\frac76};q)_{2n}(1-q^{3n+\frac{11}6})}
     {(q^{\frac{2}3};q)_{n}(q^{\frac{5}3};q)_{n}(q;q)_{2n}(q^{2};q)_{2n}(q^{\frac{8}3};q)_{2n}(1-q^{\frac{11}6})}&\\
\times\bigg\{1-
\frac{(1-q^{-n-\frac23})(1-q^{2n+\frac53})(1-q^{2n})(1-q^{2n+1})(1-q^{3n-\frac16})}
     {(1-q^{n-\frac56})(1-q^{n-\frac16})(1-q^{n+\frac{5}6})(1-q^{2n+\frac16})(1-q^{3n+\frac{11}6})}
\bigg\}&.}
\end{exam}
This provides a $q$-analogue for the following 
classical series:  
\[ \frac{32 \Gamma (\frac{1}{3})^3}{5 \pi ^2}
=\sum_{n=0}^{\infty}
\hyp{c}
{\frac{1}{6},-\frac{1}{6},\frac{5}{6},-\frac{5}{6},
\frac{11}{6},\frac{1}{12},\frac{7}{12},\frac{23}{18}}
{\\[-3mm] 1,~ 1,~ \frac{1}{2},~\frac{3}{2},~\frac{2}{3},~
\frac{4}{3},~\frac{5}{3},~\frac{5}{18}}_n
\frac{2160 n^3+1188 n^2-84 n+11}{16^n}.\]

\begin{exam}[$\boxed{a=q^{\frac32},~b=c=d=q}$ in Theorem~\ref{thm=2U}]
\pnq{\hyq{c}{q}{q,~q,~q,~q}{q^{\frac12},q^{\frac12},q^{\frac12},q^{\frac52}}_\infty
&=\sum_{n=0}^{\infty}q^{n^2+n}
\frac{(q^{\frac12};q)^2_n(q^{2};q)_n(q^{\frac32};q)_{2n}(1-q^{3n+2})}
     {(q^3;q^2)_{n}(q^{\frac52};q)^2_{2n}(-q;q)_{n}(1-q^{\frac32})}\\
&\times\bigg\{q^{n+\frac12}+
\frac{(1+q^{n+\frac12})(1-q^{2n+\frac32})^2(1-q^{3n+1})}
     {(1-q^{n+1})(1-q^{2n+\frac12})(1-q^{3n+2})}
\bigg\}.}
\end{exam}
When $q\to1$, we recover the following known series 
(cf.~\cite{kn:chu14mc}):  
\[ \frac{9 \pi ^2}{8}
=\sum_{n=0}^{\infty}
\hyp{c}
{ 1,\frac{1}{2},\frac{1}{2},\frac{1}{4},\frac{3}{4}}
{\\[-3mm] \frac{3}{2},\frac{5}{4},\frac{5}{4},\frac{7}{4},\frac{7}{4}}_n
\frac{60 n^3+111 n^2+64 n+11 }{16^n}.\]

\begin{exam}[$\boxed{a=q^{\frac12},~b=q^{\frac13},~c=q^{\frac13},~d=q^{-\frac23}}$ in Theorem~\ref{thm=2U}]
\pnq{\hyq{c}{q}{q^2,q^{\frac13},q^{\frac13},q^{-\frac23}}
               {q^{-\frac12},q^{\frac16},q^{\frac{7}6},q^{\frac{7}6}}_\infty
=\sum_{n=0}^{\infty}q^{n^2+\frac43 n}
\frac{\ff{q^{\frac13},q^{-\frac23},q^{\frac56},q^{-\frac{5}6},q^{\frac{11}6};q}{n}(q^{-\frac56};q)_{2n}(1-q^{3n-\frac23})}
     {(q^{\frac{7}6};q)_{n}(q^{-\frac12};q)_{2n}(q;q)_{2n}(q^{\frac{7}6};q)_{2n}(1-q^{n+\frac16})}&\\
\times
\bigg\{1+q^{n+\frac76}
\frac{(1-q^{n-\frac23})(1-q^{n-\frac56})(1-q^{n+\frac56})(1-q^{2n-\frac56})(1-q^{3n+\frac43})}
{(1-q^{n+\frac76})(1-q^{2n-\frac12})(1-q^{2n+1})(1-q^{2n+\frac76})(1-q^{3n-\frac23})}
     \bigg\}&.}
\end{exam}
This provides a $q$-analogue for the following 
classical series:  
\[ \frac{98 \pi ^2}{3 \Gamma (\frac{2}{3})^3}
=\sum_{n=0}^{\infty}
\hyp{c}
{\frac{1}{3},-\frac{2}{3},\frac{5}{6},-\frac{5}{6},
\frac{11}{6},\frac{10}{9},\frac{1}{12},-\frac{5}{12}}
{\\[-3mm] 1,~ \frac{3}{2},~\frac{1}{4},~\frac{3}{4},~
\frac{13}{6},~\frac{1}{9},~\frac{13}{12},~\frac{19}{12}}_n
\frac{118+45 n-1098 n^2-1080 n^3 }{16^n}.\]

\begin{exam}[$\boxed{a=q^{\frac32},~b=q^{\frac13},~c=q^{\frac13},~d=q^{\frac43}}$ in Theorem~\ref{thm=2U}]
\pnq{\hyq{c}{q}{q^2,~q^{\frac13},~q^{\frac13},~q^{\frac73}}
               {q^{\frac12},q^{\frac16},q^{\frac{13}6},q^{\frac{13}6}}_\infty
=\sum_{n=0}^{\infty}q^{n^2+\frac43 n}
\frac{\ff{q^{\frac13},q^{\frac43},q^{\frac56},q^{-\frac56},q^{\frac{11}6};q}{n}(q^{\frac16};q)_{2n}}
     {(q^{\frac{13}6};q)_{n}(q^{\frac12};q)_{2n}(q;q)_{2n}(q^{\frac{13}6};q)_{2n}}&\\
\times \frac{1-q^{3n+\frac43}}{1-q^{\frac43}}
\bigg\{1+q^{n+\frac16}
\frac{(1-q^{n+\frac43})(1-q^{n+\frac{11}6})(1-q^{2n+\frac16})}
{(1-q^{2n+\frac12})(1-q^{2n+1})(1-q^{2n+\frac{13}6})}\bigg\}&.}
\end{exam}
This provides a $q$-analogue for the following 
classical series:  
\[ \frac{637 \pi ^2}{16 \Gamma (\frac{2}{3})^3}
=\sum_{n=0}^{\infty}
\hyp{c}
{ \frac{1}{3},\frac{4}{3},\frac{5}{6},-\frac{5}{6},
\frac{11}{6},\frac{13}{9},\frac{1}{12},\frac{7}{12}}
{\\[-3mm] 1,\frac{3}{2},\frac{3}{4},~\frac{5}{4},
\frac{13}{6},\frac{4}{9},\frac{19}{12},\frac{25}{12}}_n
\frac{1080 n^3+2286 n^2+1395 n+161}{16^n}.\]

\begin{exam}[$\boxed{a=q^{\frac12},~b=q^{\frac23},~c=q^{-\frac13},~d=q^{-\frac13}}$ in Theorem~\ref{thm=2U}]
\pnq{\hyq{c}{q}{q^2,~q^{-\frac13},q^{-\frac13},q^{\frac23}}
               {q^{-\frac12},q^{-\frac16},q^{\frac{5}6},q^{\frac{11}6}}_\infty
=\sum_{n=0}^{\infty}q^{n^2+\frac23 n}
\frac{\ff{q^{\frac23},q^{-\frac13},q^{\frac{7}6},q^{-\frac76},q^{\frac{13}6};q}{n}(q^{-\frac16};q)_{2n}(1-q^{3n-\frac{1}3})}
     {(q^{\frac{5}6};q)_{n}(q^{-\frac12};q)_{2n}(q;q)_{2n}(q^{\frac{11}6};q)_{2n}(1-q^{n-\frac16})}&\\
\times\bigg\{1+q^{n+\frac56}
\frac{(1-q^{n-\frac13})(1-q^{n-\frac76})(1-q^{n+\frac{7}6})(1-q^{2n-\frac16})(1-q^{3n+\frac{5}3})}
{(1-q^{n+\frac56})(1-q^{2n-\frac12})(1-q^{2n+1})(1-q^{2n+\frac{11}6})(1-q^{3n-\frac{1}3})} \bigg\}&.}
\end{exam}
This provides a $q$-analogue for the following 
classical series:  
\[\frac{275 \pi ^2}{\Gamma (\frac{1}{3})^3}
=\sum_{n=0}^{\infty}
\hyp{c}
{\frac{2}{3},-\frac{1}{3},\frac{7}{6},-\frac{7}{6},
\frac{13}{6},\frac{11}{9},-\frac{1}{12},\frac{5}{12}}
{\\[-3mm] 1,~\frac{3}{2},~\frac{1}{4},~\frac{3}{4},~
\frac{11}{6},~\frac{2}{9},~\frac{17}{12},~\frac{23}{12}}_n
\frac{125-351 n-1602 n^2-1080 n^3}{16^n}.\]

\begin{exam}[$\boxed{a=q^{\frac32},~b=q^{\frac23},~c=q^{-\frac13},~d=q^{\frac53}}$ in Theorem~\ref{thm=2U}]
\pnq{\hyq{c}{q}{q^2,~q^{-\frac13},~q^{\frac23},~q^{\frac83}}
               {q^{\frac12},q^{-\frac16},q^{\frac{11}6},q^{\frac{17}6}}_\infty
=\sum_{n=0}^{\infty}q^{n^2+\frac23 n}
\frac{\ff{q^{\frac23},q^{\frac53},q^{\frac76},q^{-\frac76},q^{\frac{13}6};q}{n}(q^{\frac56};q)_{2n}}
     {(q^{\frac{11}6};q)_{n}(q^{\frac12};q)_{2n}(q;q)_{2n}(q^{\frac{17}6};q)_{2n}}&\\
\times \frac{1-q^{3n+\frac53}}{1-q^{\frac{5}3}}
\bigg\{1+q^{n-\frac16}
\frac{(1-q^{n+\frac53})(1-q^{n+\frac{13}6})(1-q^{2n+\frac56})}
{(1-q^{2n+\frac12})(1-q^{2n+1})(1-q^{2n+\frac{17}6})}
\bigg\}&.}
\end{exam}
This provides a $q$-analogue for the following 
classical series:  
\[ \frac{2805 \pi ^2}{4 \Gamma (\frac{1}{3})^3}
=\sum_{n=0}^{\infty}
\hyp{c}
{ \frac{2}{3},\frac{5}{3},\frac{7}{6},-\frac{7}{6},
\frac{13}{6},\frac{14}{9},\frac{5}{12},\frac{11}{12}}
{\\[-3mm] 1,\frac{3}{2},\frac{3}{4},~\frac{5}{4},
\frac{11}{6},\frac{5}{9},\frac{23}{12},\frac{29}{12}}_n
\frac{1080 n^3+2790 n^2+2151 n+478}{16^n}.\]

\begin{exam}[$\boxed{a=q^{\frac32},~b=c=q,~d=q^{\frac34}}$ in Theorem~\ref{thm=2U}]
\pnq{\hyq{c}{q}{q,~q}
               {q^{\frac12},q^{\frac32}}_\infty
=\sum_{n=0}^{\infty}q^{n^2+\frac54 n}
\frac{(q^{\frac{1}4};q)_n(q^{\frac{3}4};q)_n(q^{\frac14};q)_{2n}(1-q^{3n+\frac34})}
     {(q^{\frac32};q)_{2n}(q^{\frac54};q)_{2n}
         (-q^{\frac12};q^{\frac12})_{2n}(1-q^{n+\frac12})}&\\
\times\bigg\{1+q^{n+\frac34}
\frac{(1-q^{n+\frac14})(1-q^{n+\frac12})(1-q^{2n+\frac{1}4})(1-q^{3n+2})}
     {(1-q^{2n+1})(1-q^{2n+\frac54})(1-q^{2n+\frac32})(1-q^{3n+\frac34})}
\bigg\}&.}
\end{exam}
This provides a $q$-analogue for the following 
classical series:  
\[ 15 \pi
=\sum_{n=0}^{\infty}
\hyp{c}
{\frac{1}{2},\frac{3}{4},\frac{1}{8},\frac{5}{8}}
{\\[-3mm] \frac{3}{2},\frac{7}{4},\frac{9}{8},\frac{13}{8}}_n
\frac{120 n^2+151 n+47}{16^n}.\]

\begin{exam}[$\boxed{a=q^{\frac52},~b=q,~c=q^2,~d=q^{\frac34}}$ in Theorem~\ref{thm=2U}]
\pnq{\hyq{c}{q}{q,~q^2,~q^{\frac94}}
               {q^{\frac32},q^{\frac52},q^{\frac54}}_\infty
=\sum_{n=0}^{\infty}q^{n^2+\frac{13}4 n}
\frac{(q^{\frac12};q)_n(q^{\frac14};q)_n(q^{\frac34};q)^2_n(q^{-\frac34};q)_{2n}(1-q^{3n+\frac34})}
     {(q^{\frac52},q^{\frac74};q)_n(q^{\frac32};q)_{2n}(q^{\frac54};q)_{2n}(-q^{\frac12};q^{\frac12})_{2n}(1-q^{\frac34})}&\\
\times\bigg\{1+q^{n+\frac74}
\frac{(1-q^{n+\frac14})(1-q^{n+\frac12})(1-q^{n+\frac34})(1-q^{2n-\frac{3}4})(1-q^{3n+2})}
     {(1-q^{n+\frac74})(1-q^{2n+1})(1-q^{2n+\frac54})(1-q^{2n+\frac32})(1-q^{3n+\frac34})}
\bigg\}&.}
\end{exam}
This provides a $q$-analogue for the following 
classical series:  
\[ \frac{63 \pi }{2}
=\sum_{n=0}^{\infty}
\hyp{c}
{ \frac{1}{2},\frac{3}{4},\frac{1}{8},-\frac{3}{8}}
{\\[-3mm] \frac{5}{2},\frac{11}{4},\frac{9}{8},\frac{13}{8}}_n
\frac{120 n^2+235 n+99 }{16^n}.\]

\begin{exam}[$\boxed{a=q^{\frac32},~b=q,~c=q^{\frac56},~d=q}$ in Theorem~\ref{thm=2U}]
\pnq{\hyq{c}{q}{q,~q^2,q^{\frac56},q^{\frac76}}
               {q^{\frac12},q^{\frac32},q^{\frac43},q^{\frac53}}_\infty
=\sum_{n=0}^{\infty}q^{n^2+n}
\frac{(q;q)_n(q^{\frac13};q)^2_n(q^{\frac23};q)^2_n(q^{\frac12};q)_{2n}(1-q^{3n+1})}
     {(q^{\frac12};q)^2_n(q^{\frac32};q)_n(q^{\frac43};q)_{2n}(q^{\frac53};q)_{2n}(-q^{\frac12};q^{\frac12})_{2n}(1-q)}&\\
\times\bigg\{1+q^{n+\frac12}
\frac{(1-q^{n+\frac13})(1-q^{n+\frac23})(1-q^{n+1})(1-q^{2n+\frac{1}2})(1-q^{3n+2})}
    {(1-q^{n+\frac12})(1-q^{2n+1})(1-q^{2n+\frac43})(1-q^{2n+\frac53})(1-q^{3n+1})}
\bigg\}&.}
\end{exam}
When $q\to1$, we recover another known series 
(cf.~\cite{kn:chu14mc}):  
\[ \frac{20 \pi }{9 \sqrt{3}}
=\sum_{n=0}^{\infty}
\hyp{c}
{ 1,\frac{1}{3},\frac{2}{3},~\frac{1}{4},\frac{3}{4},\frac{8}{5}}
{\\[-3mm] \frac{3}{2},\frac{3}{2},\frac{3}{2},\frac{3}{5},\frac{7}{6},\frac{11}{6}}_n
\frac{12 n^2+15 n+4}{16^n}.\]

\begin{exam}[$\boxed{a=q^{\frac32},~b=c=q,~d=q^{\frac56}}$ in Theorem~\ref{thm=2U}]
\pnq{\hyq{c}{q}{q,~q,q^{\frac56},q^{\frac76}}
               {q^{\frac32},q^{\frac32},q^{\frac13},q^{\frac23}}_\infty
=\sum_{n=0}^{\infty}q^{n^2+\frac{7}6 n}
\frac{(q^{\frac12};q)_n(q^{\frac13};q)_n(q^{\frac23};q)_n(q^{\frac56};q)_n(1-q^{3n+\frac56})}
     {(q^{\frac32},q^{\frac23};q)_n(q^{\frac32};q)_{2n}(-q^{\frac12};q^{\frac12})_{2n}(1-q^{2n+\frac13})}&\\
\times\bigg\{1+q^{n+\frac23}
\frac{(1-q^{n+\frac13})(1-q^{n+\frac12})(1-q^{n+\frac56})(1-q^{2n+\frac{1}3})(1-q^{3n+2})}
     {(1-q^{n+\frac23})(1-q^{2n+1})(1-q^{2n+\frac43})(1-q^{2n+\frac32})(1-q^{3n+\frac56})}
\bigg\}&.}
\end{exam}
This provides a $q$-analogue for the following 
classical series:  
\[ \frac{36 \pi }{5 \sqrt{3}}
=\sum_{n=0}^{\infty}
\hyp{c}
{ \frac{1}{2},\frac{1}{3},\frac{2}{3},\frac{1}{6},\frac{5}{6},\frac{14}{9}}
{\\[-3mm] \frac{3}{2},\frac{5}{3},\frac{5}{4},\frac{7}{4},\frac{7}{6},\frac{5}{9}}_n
\frac{60 n^2+64 n+13}{16^n}.\]

\begin{exam}[$\boxed{a=q^{\frac12},~b=q^{\frac12},~c=q^{\frac12},~d=q^{\frac14}}$ in Theorem~\ref{thm=2U}]
\pnq{\hyq{c}{q}{q^{\frac12},q^{\frac32}}
               {q,~q}_\infty
=\sum_{n=0}^{\infty}q^{n^2+\frac{n}4}
\frac{(q^{\frac12};q)_n(q^{\frac14};q)_n(q^{\frac{3}4};q)_n(q^{\frac14};q)_{2n}(1-q^{3n+\frac14})}
     {(q;q)^3_{n}(q^{\frac34};q)_{2n}(-q^{\frac12};q^{\frac12})^2_{2n}(1-q^{\frac12})}&\\
\times\bigg\{1+
\frac{q^{n+\frac14}(1-q^{n+\frac14})(1-q^{n+\frac12})(1-q^{2n+\frac{1}4})(1-q^{3n+\frac{3}2})}
     {(1-q^{2n+\frac34})(1-q^{2n+1})^2(1-q^{3n+\frac14})}
\bigg\}&.}
\end{exam}
This provides a $q$-analogue for the following
remarkable series:  
\[ \frac{48}{\pi }
=\sum_{n=0}^{\infty}
\hyp{c}
{ \frac{1}{2},\frac{1}{4},\frac{3}{4},\frac{1}{8},\frac{5}{8}}
{\\[-3mm] 1,1,1,\frac{7}{8},\frac{11}{8}}_n
\frac{480 n^2+212 n+15}{16^n}.\]

\begin{exam}[$\boxed{a=q^{\frac12},~b=q^{\frac12},~c=q^{\frac12},~d=q^{\frac34}}$ in Theorem~\ref{thm=2U}]
\pnq{\hyq{c}{q}{q^{\frac12},q^{\frac12},q^{\frac{1}4},q^{\frac74}}
               {q,~q,~q^{-\frac{1}4},~q^{\frac54}}_\infty
=\sum_{n=0}^{\infty}q^{n^2-\frac{n}4}
\frac{(q^{\frac12};q)_n(q^{\frac14};q)_n(q^{\frac{3}4};q)^2_n(q^{\frac34};q)_{2n}}
     {(q;q)^3_{n}(q^{-\frac14};q)_{n}(q^{\frac{5}4};q)_{2n}(-q^{\frac12};q^{\frac12})^2_{2n}}&\\
\times\frac{1-q^{3n+\frac34}}{1-q^{\frac34}}
\bigg\{1+q^{n-\frac14}
\frac{(1-q^{n+\frac34})^2(1-q^{n+\frac12})(1-q^{2n+\frac{3}4})(1-q^{3n+\frac{3}2})}
    {(1-q^{n-\frac14})(1-q^{2n+1})^2(1-q^{2n+\frac54})(1-q^{3n+\frac34})}
\bigg\}&.}
\end{exam}
This provides a $q$-analogue for the following 
classical series:  
\[ \frac{80}{3 \pi }
=\sum_{n=0}^{\infty}
\hyp{c}
{ \frac{1}{2},\frac{1}{4},\frac{3}{4},\frac{3}{8},\frac{7}{8}}
{\\[-3mm] 1,1,1,\frac{9}{8},\frac{13}{8}}_n
\frac{640 n^3+560 n^2+112 n+7}{16^n}.\]

\begin{exam}[$\boxed{a=q^{\frac12},~b=q^{\frac12},~c=q^{\frac12},~d=q^{\frac13}}$ in Theorem~\ref{thm=2U}]
\pnq{\hyq{c}{q}{q^{\frac12},q^{\frac12},q^{\frac{2}3},q^{\frac43}}
               {q,~q,~q^{\frac16},~q^{\frac{5}6}}_\infty
=\sum_{n=0}^{\infty}q^{n^2+\frac{n}6}
\frac{(q^{\frac12};q)_n(q^{\frac13};q)^2_n(q^{\frac{2}3};q)_n(q^{\frac13};q)_{2n}}
     {(q;q)^3_{n}(q^{\frac{1}6};q)_{n}(q^{\frac{5}6};q)_{2n}(-q^{\frac12};q^{\frac12})^2_{2n}}&\\
\times\frac{1-q^{3n+\frac13}}{1-q^{\frac13}}
\bigg\{1+
\frac{q^{n+\frac16}(1-q^{n+\frac13})^2(1-q^{n+\frac12})(1-q^{2n+\frac{1}3})(1-q^{3n+\frac{3}2})}
     {(1-q^{n+\frac16})(1-q^{2n+\frac56})(1-q^{2n+1})^2(1-q^{3n+\frac13})}
\bigg\}&.}
\end{exam}
When $q\to1$, we recover the following known series 
(cf.~\cite{kn:chu14mc}):  
\[ \frac{15 \sqrt{3}}{\pi }
=\sum_{n=0}^{\infty}
\hyp{c}
{ \frac{1}{2},~\frac{1}{3},\frac{1}{3},\frac{2}{3},\frac{2}{3}}
{\\[-3mm] 1,1,1,\frac{11}{12},\frac{17}{12}}_n
\frac{135 n^2+75 n+8}{16^n}.\]

\begin{exam}[$\boxed{a=q^{\frac12},~b=q^{\frac12},~c=q^{\frac16},~d=q^{\frac12}}$ in Theorem~\ref{thm=2U}]
\pnq{\hyq{c}{q}{q^{\frac32},q^{\frac52},q^{\frac{7}6},q^{\frac{11}6}}
               {q,~q^{2},~q^{\frac{5}3},~q^{\frac73}}_\infty
=\sum_{n=0}^{\infty}q^{n^2}
\frac{\ff{q^{\frac12},q^{\frac16},q^{\frac56},q^{\frac76},q^{\frac{11}6};q}{n}(q^{\frac32};q)_{2n}(1-q^{3n+\frac{3}2})}
     {(q;q)^3_{n}(q^{\frac{5}3};q)_{2n}(q^{\frac{7}3};q)_{2n}(-q;q^{\frac12})_{2n}(1-q^{\frac32})}&\\
\times\bigg\{q^{n}+
\frac{(1-q^{n})(1-q^{2n+\frac23})(1-q^{2n+1})(1-q^{2n+\frac43})(1-q^{3n+\frac12})}
     {(1-q^{n+\frac16})(1-q^{n+\frac12})(1-q^{n+\frac56})(1-q^{2n+\frac{1}2})(1-q^{3n+\frac{3}2})}
\bigg\}&.}
\end{exam}
This provides a $q$-analogue for the following 
classical series:  
\[ \frac{256}{3 \pi \sqrt{3}}
=\sum_{n=0}^{\infty}
\hyp{c}
{ \frac{1}{2},\frac{1}{4},\frac{3}{4},\frac{1}{6},\frac{5}{6}}
{\\[-3mm] 1,1,1,\frac{4}{3},\frac{5}{3} }_n
\frac{720 n^3+804 n^2+236 n+15}{16^n}.\]

\begin{exam}[$\boxed{a=q^{\frac12},~b=q^{\frac12},~c=q^{\frac12},~d=q^{\frac56}}$ in Theorem~\ref{thm=2U}]
\pnq{\hyq{c}{q}{q^{\frac12},q^{\frac12},q^{\frac16},q^{\frac{11}6}}
               {q,~q,q^{-\frac13},q^{\frac43}}_\infty
=\sum_{n=0}^{\infty}q^{n^2-\frac{n}3}
\frac{(q^{\frac12};q)_n(q^{\frac16};q)_n(q^{\frac56};q)^2_n(q^{\frac56};q)_{2n}}
     {(q;q)^3_{n}(q^{-\frac{1}3};q)_{n}(q^{\frac43};q)_{2n}(-q^{\frac12};q^{\frac12})^2_{2n}}&\\
\times\frac{1-q^{3n+\frac56}}
     {1-q^{\frac56}}
\bigg\{1+q^{n-\frac13}
\frac{(1-q^{n+\frac12})(1-q^{n+\frac56})^2(1-q^{2n+\frac{5}6})(1-q^{3n+\frac{3}2})}
     {(1-q^{n-\frac13})(1-q^{2n+1})^2(1-q^{2n+\frac43})(1-q^{3n+\frac56})}
\bigg\}&.}
\end{exam}
When $q\to1$, we recover the following known series 
(cf.~\cite{kn:chu14mc}):  
\[ \frac{384}{\pi\sqrt{3} }
=\sum_{n=0}^{\infty}
\hyp{c}
{ \frac{1}{2},\frac{5}{6},\frac{5}{6},\frac{5}{12},\frac{11}{12}}
{\\[-3mm] 1,1,1,~ \frac{2}{3},\frac{5}{3}}_n
\frac{1080 n^2+798 n+55}{16^n}.\]

\subsection{} \
Analogously, letting $\boxed{\del=1}$ in Theorem~\ref{thm=V}, we obtain,
under the partition pattern \eqref{p2pattern}, another summation formula.
\begin{thm}[Nonterminating series identity]\label{thm=2V}
\pnq{\hyq{c}{q}{b,~c,~qd,~qa^2/bcd}{qa/b,qa/c,a/d,bcd/a}_\infty
=1+\sum_{n=1}^{\infty}W_n(a,b,c,d)
\Big(\frac{q^na^2}{bd}\Big)^n&\\
\times\frac{\ff{bd/a;q}{2n}\ff{b,d,qa/bc,bc/a,qa/cd,cd/a;q}{n}}
{\ff{q,qa/c,bcd/a;q}{2n}\ff{qa/b,a/d;q}{n}},&}
where the weight function is given by
\pnq{&W_n(a,b,c,d)=\frac{1-q^{3n}d}{1-d}\\
&\times\bigg\{1-
\frac{(1-q^{2n})(1-q^{-n}b/a)(1-q^{2n}a/c)(1-q^{2n-1}bcd/a)(1-q^{3n-2}b)}
{(1-q^na/cd)(1-q^{n-1}b)(1-q^{n-1}bc/a)(1-q^{2n-1}bd/a)(1-q^{3n}d)}\bigg\}.}
\end{thm}

By specifying the parameters $\{a,b,c,d\}$ concretely, we can deduce
from Theorem~\ref{thm=2V} a number of $q$-series identities that are
$q$-analogues of classical series counterparts with the convergence
rate ``$\frac1{16}$". Twelve formulae are recorded below as examples.

\begin{exam}[$\boxed{a=b=c=d=q^{\frac12}}$ in Theorem~\ref{thm=2V}]
\pnq{\hyq{c}{q}{q^{\frac12},q^{\frac12},q^{\frac12},q^{\frac52}}{q,~q,~q,~q}_\infty
&=\sum_{n=1}^{\infty}q^{n^2}
\frac{(q^{\frac12};q)^3_n(q^{\frac12};q)_{2n}}
     {(q;q)^5_n(-q^{\frac12};q^{\frac12})^3_{2n}}
\frac{(1-q^n)(1-q^{3n+\frac12})}
     {(1-q^{1/2})(1-q^{3/2})}\\
&\times\bigg\{1-
\frac{(1-q^{-n})(1-q^{2n})^3(1-q^{3n-\frac32})}
     {(1-q^{2n-\frac12})(1-q^{n-\frac12})^3(1-q^{3n+\frac12})}\bigg\}.}
\end{exam}
This provides a $q$-analogue for the following 
classical series:  
\[\frac{2}{\pi^2}
=\sum_{n=1}^{\infty}\hyp{c}
{\frac12,-\frac12,-\frac12,\frac14,-\frac14}
{\\[-3mm]1,~1,~1,~1,~1}_n
\frac{n(1-6n)(1+4n-4n^2+80n^3)}{16^n}.\]

\begin{exam}[$\boxed{a=q^{\frac12},~b=q^{\frac56},~c=q^{\frac56},~d=q^{\frac56}}$ in Theorem~\ref{thm=2V}]
\pnq{\hyq{c}{q}{q^{\frac12},q^{\frac56},~q^{\frac56},~q^{\frac{5}6}}
               {q^{2},q^{-\frac13},q^{\frac23},q^{\frac{2}3}}_\infty
=\sum_{n=0}^{\infty}q^{n^2-\frac{2}3n}
\frac{(q^{-\frac16};q)^2_n(q^{\frac56};q)^2_n(q^{\frac{7}6};q)^2_n(q^{\frac76};q)_{2n}}
     {(q^{-\frac{1}3};q)_{n}(q^{\frac{2}3};q)_{n}(q;q)_{2n}(q^{2};q)_{2n}(q^{\frac{2}3};q)_{2n}}&\\
\times\frac{1-q^{3n+\frac{5}6}}
     {1-q^{-\frac12}}
\bigg\{1-
\frac{(1-q^{\frac13-n})(1-q^{2n-\frac13})(1-q^{2n})(1-q^{2n+1})(1-q^{3n-\frac76})}
     {(1-q^{n-\frac76})(1-q^{n-\frac16})(1-q^{n+\frac{1}6})(1-q^{2n+\frac16})(1-q^{3n+\frac{5}6})}
\bigg\}&.}
\end{exam}
This provides a $q$-analogue for the following 
classical series:  
\[ \frac{63 \Gamma (\frac{1}{3})^3}{10 \pi ^2}
=\sum_{n=0}^{\infty}
\hyp{c}
{ \frac{1}{6},-\frac{1}{6},\frac{5}{6},\frac{7}{6},
-\frac{7}{6},\frac{1}{12},\frac{7}{12},\frac{17}{18}}
{\\[-3mm] 1,~ 1,~ \frac{1}{2},\frac{3}{2},-\frac{1}{3},
\frac{1}{3},\frac{2}{3},-\frac{1}{18}}_n
\frac{432 n^3-108 n^2-132 n+7}{16^n}.\]

\begin{exam}[$\boxed{a=q^{\frac52},~b=q,~c=q^2,~d=q}$ in Theorem~\ref{thm=2V}]
\pnq{\hyq{c}{q}{q,~q,~q^2,~q^2}{q^{\frac12},q^{\frac32},q^{\frac32},q^{\frac52}}_\infty
=1+q^{\frac12}+\sum_{n=1}^{\infty}q^{n^2+3n}
\frac{(q^{\frac12};q)^2_n(q;q)_n(q^{-\frac12};q)_{2n}}
     {(q^{\frac52};q)_{n}(q^{\frac32};q)^2_{2n}(-q^{\frac12};q^{\frac12})_{2n}}&\\
\times\frac{1-q^{3n+1}}{1-q^{n+\frac12}}
\bigg\{1-
\frac{(1+q^{n})(1-q^{-n-\frac32})(1-q^{2n+\frac12})^2(1-q^{3n-1})}
     {(1-q^{n-\frac12})^2(1-q^{2n-\frac32})(1-q^{3n+1})}
\bigg\}&.}
\end{exam}
This provides a $q$-analogue for the following 
classical series:  
\[\frac{3 \pi ^2}{32}
=1+\sum_{n=1}^{\infty}
\hyp{c}
{1,-\frac{1}{2},-\frac{1}{2},-\frac{1}{4},-\frac{3}{4}}
{\\[-3mm]\frac{5}{2},~\frac{3}{4},~\frac{3}{4},~\frac{5}{4},~\frac{5}{4}}_n
\frac{3+3n-22 n^2-40 n^3}{16^n}.\]

\begin{exam}[$\boxed{a=q^{\frac72},~b=q,~c=q^2,~d=q}$ in Theorem~\ref{thm=2V}]
\pnq{\hyq{c}{q}{q,~q^2,~q^2,~q^4}{q^{\frac12},q^{\frac52},q^{\frac52},q^{\frac72}}_\infty
=1+\sum_{n=1}^{\infty}q^{n^2+5n}
\frac{(q^{-\frac12};q)^2_n(q;q)_n(q^{-\frac32};q)_{2n}(1-q^{3n+1})}
     {(q^{\frac72};q)_{n}(q^{\frac12};q)_{2n}(q^{\frac52};q)_{2n}(-q^{\frac12};q^{\frac12})_{2n}(1-q)}&\\
\times\frac{(1-q^{\frac32})(1-q^{n+\frac12})} {(1-q^{\frac12})(1-q^{n+\frac32})}
\bigg\{1-
\frac{(1+q^{n})(1-q^{-n-\frac52})(1-q^{2n-\frac12})(1-q^{2n+\frac32})(1-q^{3n-1})}
     {(1-q^{n-\frac32})(1-q^{n+\frac12})(1-q^{2n-\frac52})(1-q^{3n+1})}
\bigg\}&.}
\end{exam}
This provides a $q$-analogue for the following 
classical series:  
\[\frac{15 \pi ^2}{256}
=\frac13+\sum_{n=1}^{\infty}
\hyp{c}
{1,-\frac{1}{2},-\frac{3}{2},-\frac{3}{4},-\frac{5}{4}}
{\\[-3mm] \frac{7}{2},~ \frac{1}{4},~ \frac{3}{4},~ \frac{5}{4},~ \frac{7}{4}}_n
\frac{8 n^3+2 n^2-3 n+1}{16^n}.\]

\begin{exam}[$\boxed{a=q^{\frac52},~b=q^{2},~c=q^{2},~d=q^{\frac34}}$ in Theorem~\ref{thm=2V}]
\pnq{\hyq{c}{q}{q,q^{2}}{q^{\frac32},q^{\frac32}}_\infty
=\sum_{n=0}^{\infty}q^{n^2+\frac94 n}
\frac{(q^{-\frac12};q)_n(q^{\frac14};q)_n(q^{\frac34};q)_n(q^{\frac14};q)_{2n}(1-q^{n+1})}
     {(q^{\frac12};q)_{n}(q^{\frac{3}2};q)_{2n}(q^{\frac94};q)_{2n}(-q^{\frac12};q^{\frac12})_{2n}(1-q^{n+\frac34})}&\\
\times\frac{1-q^{3n+\frac{3}4}}
     {1-q^{\frac54}}
\bigg\{1+
\frac{q^{-n-\frac12}(1-q^{2n})(1-q^{3n})(1-q^{2n+\frac12})(1-q^{2n+\frac54})}
     {(1-q^{n-\frac14})(1-q^{n+1})(1-q^{2n-\frac34})(1-q^{3n+\frac{3}4})}
\bigg\}&.}
\end{exam}
This provides a $q$-analogue for the following 
classical series:  
\[ \frac{5 \pi }{16}
=\sum_{n=0}^{\infty}
\hyp{c}
{ -\frac{1}{2},-\frac{1}{4},\frac{1}{8},-\frac{3}{8} }
{\\[-3mm] \frac{1}{2},~\frac{3}{4},~\frac{9}{8},~\frac{13}{8} }_n
\frac{40 n^2-7 n+1}{16^n}.\]

\begin{exam}[$\boxed{a=q^{\frac52},~b=q^{2},~c=q^{2},~d=q^{\frac14}}$ in Theorem~\ref{thm=2V}]
\pnq{\hyq{c}{q}{q,~q^{2}}{q^{\frac32},q^{\frac32}}_\infty
=\sum_{n=0}^{\infty}q^{n^2+\frac{11}4 n}
\frac{(q^{-\frac12};q)_n(q^{\frac14};q)_n(q^{-\frac14};q)_n(q^{-\frac14};q)_{2n}(1-q^{3n+\frac{1}4})}
     {(q^{\frac12};q)_{n}(q^{\frac{3}2};q)_{2n}(q^{\frac74};q)_{2n}(-q^{\frac12};q^{\frac12})_{2n}(1-q^{\frac14})}&\\
\times
\frac{1-q^{n+1}}{1-q^{n+\frac54}}
\bigg\{1+
\frac{q^{-n-\frac12}(1-q^{2n})(1-q^{3n})(1-q^{2n+\frac12})(1-q^{2n+\frac34})}
     {(1-q^{n+\frac14})(1-q^{n+1})(1-q^{2n-\frac54})(1-q^{3n+\frac{1}4})}
\bigg\}&.}
\end{exam}
This provides a $q$-analogue for the following 
classical series:  
\[ \frac{25 \pi }{16}
=\sum_{n=0}^{\infty}
\hyp{c}
{ -\frac{1}{2},\frac{1}{4},-\frac{1}{8},-\frac{5}{8}}
{\\[-3mm] \frac{1}{2},~\frac{9}{4},~\frac{7}{8},~\frac{11}{8}}_n
\frac{120 n^2+77 n+5}{16^n}.\]

\begin{exam}[$\boxed{a=q^{\frac32},~b=q,~c=q,~d=q^{\frac14}}$ in Theorem~\ref{thm=2V}]
\pnq{\hyq{c}{q}{q,~q,~q^{\frac74}}
               {q^{\frac32},q^{\frac32},q^{\frac34}}_\infty
=1+\sum_{n=1}^{\infty}q^{n^2+\frac74 n}
\frac{(q^{\frac12};q)_n(q^{\frac14};q)_n(q^{-\frac14};q)_n(q^{-\frac14};q)_{2n}}
     {(q^{\frac32};q)_{n}(q^{\frac{3}2};q)_{2n}
         (q^{\frac34};q)_{2n}(-q^{\frac12};q^{\frac12})_{2n}}&\\
\times\frac{1-q^{3n+\frac14}}{1-q^{\frac14}}
\bigg\{1-
\frac{(1+q^{n})(1-q^{-n-\frac12})(1-q^{2n-\frac14})(1-q^{2n+\frac12})(1-q^{3n-1})}
     {(1-q^{n-\frac12})(1-q^{n+\frac14})(1-q^{2n-\frac54})(1-q^{3n+\frac14})}
\bigg\}&.}
\end{exam}
This provides a $q$-analogue for the following 
classical series:  
\[\frac{5 \pi }{9}
=\frac53+\sum_{n=1}^{\infty}
\hyp{c}
{ -\frac{1}{2},\frac{1}{4},-\frac{1}{8},-\frac{5}{8}}
{\\[-3mm] \frac{3}{2},~\frac{5}{4},~\frac{3}{8},~\frac{7}{8}}_n
\frac{7-6n-80 n^2}{16^n}.\]

\begin{exam}[$\boxed{a=q^{\frac52},~b=q,~c=q^2,~d=q^{\frac14}}$
        in Theorem~\ref{thm=2V}]
\pnq{\hyq{c}{q}{q,~q^2,~q^{\frac54},q^{\frac{11}4}}
               {q^{\frac32},q^{\frac52},q^{\frac34},q^{\frac94}}_\infty
=1+\sum_{n=1}^{\infty}q^{n^2+\frac{15}4 n}
\frac{(q^{\frac12};q)_n(q^{\frac14};q)_n(q^{-\frac14};q)_n(q^{\frac54};q)_n(q^{-\frac54};q)_{2n}}
     {(q^{\frac52};q)_{n}(q^{\frac94};q)_n(q^{\frac{3}2};q)_{2n}
         (q^{\frac34};q)_{2n}(-q^{\frac12};q^{\frac12})_{2n}}&\\
\times\frac{1-q^{3n+\frac14}}{1-q^{\frac14}}
\bigg\{1-
\frac{(1+q^{n})(1-q^{-n-\frac32})(1-q^{2n-\frac14})(1-q^{2n+\frac12})(1-q^{3n-1})}
     {(1-q^{n-\frac12})(1-q^{n+\frac14})(1-q^{2n-\frac94})(1-q^{3n+\frac14})}
\bigg\}&.}
\end{exam}
This provides a $q$-analogue for the following 
classical series:  
\[\frac{15 \pi }{14}
=3+\sum_{n=1}^{\infty}
\hyp{c}
{ -\frac{1}{2},\frac{1}{4},-\frac{5}{8},-\frac{9}{8}}
{\\[-3mm] \frac{5}{2},~\frac{9}{4},~\frac{3}{8},~\frac{7}{8} }_n
\frac{19-62 n-80 n^2}{16^n}.\]

\begin{exam}[$\boxed{a=q^{\frac12},~b=q,~c=q^{\frac16},~d=q}$ in Theorem~\ref{thm=2V}]
\pnq{\hyq{c}{q}{q,~q^2,~q^{\frac16},q^{-\frac16}}
               {q^{-\frac12},q^{\frac12},q^{\frac43},q^{\frac53}}_\infty
=1+
\sum_{n=1}^{\infty}q^{n^2-n}
\frac{(q;q)_n(q^{\frac13};q)^2_n(q^{\frac23};q)^2_n(q^{\frac32};q)_{2n}}
     {(q^{\frac12};q)^2_{n}(q^{-\frac12};q)_{n}(q^{\frac{4}3};q)_{2n}
         (q^{\frac53};q)_{2n}(-q^{\frac12};q^{\frac12})_{2n}}&\\
\times\frac{1-q^{3n+1}}{1-q}
\bigg\{1-
\frac{(1+q^{n})(1-q^{\frac12-n})(1-q^{2n+\frac13})(1-q^{2n+\frac23})(1-q^{3n-1})}
     {(1-q^{n-\frac23})(1-q^{n-\frac13})(1-q^{2n+\frac12})(1-q^{3n+1})}
\bigg\}&.}
\end{exam}
This provides a $q$-analogue for the following 
classical series:  
\[\frac{4 \pi }{81 \sqrt{3}}
=\frac23+\sum_{n=1}^{\infty}
\hyp{c}
{ 1,~ \frac{2}{3},-\frac{2}{3},\frac{1}{4},\frac{3}{4}}
{\\[-3mm]\frac{1}{2},\frac{1}{2},-\frac{1}{2},\frac{5}{6},\frac{7}{6}}_n
\frac{2+7n-20 n^2}{16^n}.\]

\begin{exam}[$\boxed{a=q^{\frac12},~b=q,~c=q^{\frac16},~d=q^2}$ in Theorem~\ref{thm=2V}]
\pnq{\hyq{c}{q}{q,~q^3,~q^{\frac16},q^{-\frac76}}
               {q^{-\frac32},q^{\frac12},q^{\frac43},q^{\frac83}}_\infty
=1+
\sum_{n=1}^{\infty}q^{n^2-2n}
\frac{(q^2;q)_n(q^{\frac13};q)_n(q^{\frac23};q)_n(q^{-\frac23};q)_n(q^{\frac53};q)_n(q^{\frac52};q)_{2n}}
     {(q^{\frac12};q)^2_{n}(q^{-\frac32};q)_{n}(q^{\frac{4}3};q)_{2n}
         (q^{\frac83};q)_{2n}(-q^{\frac12};q^{\frac12})_{2n}}&\\
\times\frac{1-q^{3n+2}}
     {1-q^2}
\bigg\{1-
\frac{(1+q^{n})(1-q^{\frac12-n})(1-q^{2n+\frac13})(1-q^{2n+\frac53})(1-q^{3n-1})}
     {(1-q^{n-\frac53})(1-q^{n-\frac13})(1-q^{2n+\frac32})(1-q^{3n+2})}
\bigg\}&.}
\end{exam}
This provides a $q$-analogue for the following 
classical series:  
\[\frac{350 \pi }{243 \sqrt{3}}
=30+\sum_{n=1}^{\infty}
\hyp{c}
{2,~\frac{5}{3},-\frac{5}{3},\frac{3}{4},\frac{5}{4}}
{\\[-3mm]\frac{1}{2},\frac{1}{2},-\frac{3}{2},\frac{7}{6},\frac{11}{6}}_n
\frac{40-n-60 n^2}{16^n}.\]

\begin{exam}[$\boxed{a=q^{\frac12},~b=q^{\frac12},~c=q^{\frac16},~d=q^{\frac12}}$ in Theorem~\ref{thm=2V}]
\pnq{\hyq{c}{q}{q^{\frac12},q^{\frac32},q^{\frac56},q^{\frac{7}6}}
               {q,~q,~q^{\frac23},q^{\frac43}}_\infty
=\sum_{n=0}^{\infty}q^{n^2}
\frac{(q^{\frac12};q)_n(q^{\frac16};q)^2_n(q^{\frac56};q)^2_n(q^{\frac12};q)_{2n}(1-q^{3n+\frac{1}2})}
     {(q;q)^3_{n}(q^{\frac{2}3};q)_{2n}(q^{\frac43};q)_{2n}(-q^{\frac12};q^{\frac12})_{2n}(1-q^{\frac16})}&\\
\times\frac{1-q^n}{1-q^{\frac12}}
\bigg\{1-
\frac{(1-q^{-n})(1-q^{2n})(1-q^{2n-\frac13})(1-q^{2n+\frac13})(1-q^{3n-\frac32})}
     {(1-q^{n-\frac12})(1-q^{n-\frac16})(1-q^{n-\frac56})(1-q^{2n-\frac12})(1-q^{3n+\frac{1}2})}
\bigg\}&.}
\end{exam}
This provides a $q$-analogue for the following 
classical series:  
\[ \frac{10\sqrt{3}}{9  \pi }
=\sum_{n=0}^{\infty}
\hyp{c}
{ \frac{1}{2},\frac{1}{4},-\frac{1}{4},\frac{5}{6},-\frac{5}{6}}
{\\[-3mm] 1,~ 1,~ 1,~ \frac{1}{3},~\frac{2}{3}}_n
\frac{n (120 n^2-26 n+5)}{16^n}.\]

\begin{exam}[$\boxed{a=q^{\frac12},~b=q^{\frac12},~c=q^{\frac12},~d=q^{-\frac13}}$ in Theorem~\ref{thm=2V}]
\pnq{\hyq{c}{q}{q^{\frac12},q^{\frac12},q^{\frac23},q^{\frac{4}3}}
               {q,~q,~q^{\frac16},~q^{\frac56}}_\infty
=\sum_{n=0}^{\infty}q^{n^2+\frac{5}6n}
\frac{(q^{\frac12};q)_n(q^{-\frac13};q)^2_n(q^{\frac{4}3};q)_n(q^{-\frac13};q)_{2n}}
     {(q;q)^3_{n}(q^{\frac{5}6};q)_{n}(q^{\frac16};q)_{2n}(-q^{\frac12};q^{\frac12})^2_{2n}}&\\
\times\frac{1-q^{3n-\frac{1}3}}
     {1-q^{-\frac13}}
\bigg\{1-
\frac{(1-q^{-n})(1-q^{2n})^2(1-q^{2n-\frac56})(1-q^{3n-\frac32})}
     {(1-q^{n-\frac12})^2(1-q^{n+\frac13})(1-q^{2n-\frac43})(1-q^{3n-\frac{1}3})}
\bigg\}&.}
\end{exam}
This provides a $q$-analogue for the following 
classical series:  
\[ \frac{6 \sqrt{3}}{\pi }
=\sum_{n=0}^{\infty}
\hyp{c}
{ -\frac{1}{2},\frac{1}{3},-\frac{1}{3},-\frac{1}{3},-\frac{2}{3}}
{\\[-3mm] 1,~ 1,~ 1,~ \frac{1}{12},~\frac{7}{12}}_n
\frac{135 n^3-48 n^2-7 n+2}{16^n}.\]


\section{Series from Triplicate Inversions}
For all $n\in\mb{N}_0$, it is well known that
$\boxed{n=\pavv{\tfrac{n}3}+\pavv{\tfrac{1+n}3}+\pavv{\tfrac{2+n}3}}$.
Then it is not difficult to check that Jackson's formula \eqref{jackson}
is equivalent to the following one
\pnq{\Ome_n\big(a;q^{\pav{\tfrac{n}3}}b,q^{\pav{\tfrac{1+n}3}}c,q^{\pav{\tfrac{2+n}3}}d\big)
=\hyq{c}{q}{qa/cd}{b/a}_{\pavv{\frac{n}3}}
\hyq{c}{q}{qa/bd}{c/a}_{\pavv{\frac{1+n}3}}
\hyq{c}{q}{qa/bc}{d/a}_{\pavv{\frac{2+n}3}}&\\
\times\hyq{c}{q}{qa}{bcd/a}_n
\hyq{c}{q}{bc/a}{qa/d}_{\pavv{\frac{2n}3}}
\hyq{c}{q}{bd/a}{qa/c}_{\pavv{\frac{1+2n}3}}
\hyq{c}{q}{cd/a}{qa/b}_{\pavv{\frac{2+2n}3}}&}
with its parameters subject to $\boxed{qa^2=bcde}$.
The partition pattern corresponding to \eqref{p-pattern} becomes
\pq{\label{p3pattern}\boxed{\Lam=3:\mult{cccc}
{\veps_b=0,&\veps_c=1,&\veps_d=2,&\veps_e=0;\\
\lam_b=1,&\lam_c=1,&\lam_d=1,&\lam_e=0.}}}

\subsection{} \
Specifying with $\boxed{\del=0}$ and the above partition pattern in Theorem~\ref{thm=U},
we can express the corresponding formula as in the following theorem.
\begin{thm}[Nonterminating series identity]\label{thm=3U}
\pnq{\hyq{c}{q}{b,~qc,~qd,~qa^2/bcd}{qa/b,qa/c,qa/d,bcd/a}_\infty
=\sum_{n=0}^{\infty}W_n(a,b,c,d)
\Big(\frac{q^{3n}a^3}{bcd}\Big)^n&\\
\times\frac{\ff{b,c,qd,qa/bc,qa/bd,qa/cd;q}{n}\ff{bc/a,bd/a,cd/a;q}{2n}}
{\ff{qa/b,qa/c,qa/d;q}{2n}\ff{q,bcd/a;q}{3n}},&}
where the weight function is given by
\pnq{W_n&(a,b,c,d)
=\frac{(1-q^{2n}bd/a)(1-q^{1+n}a/(bc))(1-q^{2n}cd/a)(1-q^{1+4n}c)}
      {(1-q^{1+3n})(1-c)(1-q^{1+2n}a/b)(1-q^{3n}bcd/a) }\\
\times~&\bigg\{q^{2n}\frac{a}{d}+
\frac{(1-q^{2n}a/d)(1-q^{3n}bcd/a)(1-q^{4n}d)(1-q^{1+2n}a/b)(1-q^{1+3n})}
     {(1-q^{n}d) (1-q^{2n}bd/a) (1-q^{2n}cd/a) (1-q^{1+n}a/(bc))(1-q^{1+4n}c)}&\\
+&~\frac{q^{1+4n}a^2(1-q^{n}c) (1-q^{2n}bc/a) (1-q^{1+n}a/(bd)) (1-q^{1+2n}cd/a) (1-q^{2+4n}b)}
     {cd(1-q^{1+2n}a/c) (1-q^{1+2n}a/d)(1-q^{1+3n}bcd/a) (1-q^{2+3n})(1-q^{1+4n}c)}\bigg\}.}
\end{thm}

According to this theorem, by assigning concrete values for the quadruplet
$\{a,b,c,d\}$, we can derive numerous $q$-series identities after some
routine simplifications. Even though the $q$-series so obtained are quite
complicated in general, they can dramatically be reduced, in some cases,
to simple ones by means of the ``reverse bisection method". We take
the formulae from Examples~\ref{v1x3a} and~\ref{v3x1a} to illustrate
how this approach works.

When $\boxed{a=q^{\frac43},~b=c=q,~d=q^{\frac56}}$, the corresponding
formula in Theorem~\ref{thm=3U} reads as
\pnq{\hyq{c}{q}{q,~q,~q^{\frac56},~q^{\frac56}}
           {q^{\frac32},q^{\frac32},q^{\frac43},q^{\frac43}}_\infty
&=\sum_{n=0}^{\infty}q^{3n^2+\frac76 n}
\frac{\ff{q,q,q^{\frac12},q^{\frac12},q^{\frac13},q^{\frac56};q}{n}
     \ff{q^{\frac12},q^{\frac12},q^{\frac23};q}{2n}}
{\ff{q^{\frac32},q^{\frac43},q^{\frac43};q}{2n}\ff{q,q^{\frac32};q}{3n}}
\md{W}_n(q^{\frac43},q,q,q^{\frac56}),}
where $\md{W}_n$ is a rational function of $q^n$ given explicitly by
\pnq{\md{W}_n(q^{\frac43},q,q,q^{\frac56})
&=\frac{(1-q^{n+\frac13})(1-q^{n+\frac56})(1-q^{2n+\frac12})^2(1-q^{4n+2})}
{(1-q^{2n+\frac43})(1-q^{3n+1})(1-q^{3n+\frac32})}\\
&\times\bigg\{q^{2n+\frac12}
+\frac{(1-q^{2n+\frac43})(1-q^{3n+1})(1-q^{3n+\frac32})(1-q^{4n+\frac56})}
      {(1-q^{n+\frac13})(1-q^{n+\frac56})(1-q^{2n+\frac12})(1-q^{4n+2})}\\
&+\frac{q^{4n+\frac{11}6}(1-q^{n+\frac12})(1-q^{n+1})(1-q^{2n+\frac23})(1-q^{4n+3})}
     {(1-q^{2n+\frac43})(1-q^{3n+2})(1-q^{3n+\frac52})(1-q^{4n+2})}\bigg\}.}

By factorization with \emph{Mathematica} commands, we can check that
\[\md{P}(q^n):=\frac{(1-q^{2n+\frac43})^2(1-q^{3n+1})(1-q^{3n+2})(1-q^{3n+\frac32})(1-q^{3n+\frac52})}
{(1-q^{n+\frac13})(1-q^{2n+\frac12})}
\md{W}_n(q^{\frac43},q,q,q^{\frac56})\]
is a polynomial of degree 19 in $q^n$. Then the preceding formula
in question can be rewritten, by substitution, as
\pq{\label{v1x3PP}\pp{c}{
\hyq{c}{q}{q,~q,~q^{\frac56},q^{\frac56}}
         {q^{\frac12},q^{\frac32},q^{\frac13},q^{\frac43}}_\infty
&=\sum_{n=0}^{\infty}q^{3n^2+\frac{7}6 n}
\frac{\ff{q^{\frac12},q^{\frac23};q}{2n}
\ff{q^{\frac12},q^{\frac12},q^{\frac56};q^{\frac12}}{2n}}
{(1-q)\ff{q^{\frac43},q^{\frac43};q}{2n+1}
\ff{q^{\frac32};q^{\frac12}}{6n+3}}\md{P}(q^n).}}
Supposing that $\md{Q}(y)$ is a polynomial (to be determined),
the series on the right suggests us to introduce the hypergeometric term
\[\md{T}_n:=q^{\frac34n^2+\frac{7}{12}n}
\frac{\ff{q^{\frac12},q^{\frac23};q}{n}
\ff{q^{\frac12},q^{\frac12},q^{\frac56};q^{\frac12}}{n}}
{(1-q)\ff{q^{\frac43},q^{\frac43};q}{n}
\ff{q^{\frac32};q^{\frac12}}{3n}}\md{Q}(q^{\frac{n}2}).\]
By computing the difference of two consecutive terms
\pnq{\md{T}_{2n}\pm \md{T}_{2n+1}
&=q^{3n^2+\frac{7}6 n}
\frac{\ff{q^{\frac12},q^{\frac23};q}{2n}
\ff{q^{\frac12},q^{\frac12},q^{\frac56};q^{\frac12}}{2n}}
{(1-q)\ff{q^{\frac43},q^{\frac43};q}{2n+1}
\ff{q^{\frac32};q^{\frac12}}{6n+3}}\\
&\times\Big\{\md{Q}(q^{n})(1-q^{\frac43+2n})^2(1-q^{\frac32+3n})(1-q^{2+3n})(1-q^{\frac52+3n})\\
&\pm q^{3n+\frac43}\md{Q}(q^{n+\frac12})
     (1-q^{\frac12+n})^2(1-q^{\frac56+n})(1-q^{\frac12+2n})(1-q^{\frac23+2n})\Big\},}
we can formally consider the series in \eqref{v1x3PP} as the following ``bisection series"
\pq{\label{v1x3TT}
\hyq{c}{q}{q,~q,~q^{\frac56},~q^{\frac56}}
    {q^{\frac12},q^{\frac32},q^{\frac13},q^{\frac43}}_\infty
=\sum_{n=0}^{\infty}(\pm1)^n\md{T}_n
=\sum_{n=0}^{\infty}\Big\{\md{T}_{2n}\pm \md{T}_{2n+1}\Big\},}
where the signs `$\pm$' indicate that the series can be positive
or alternating. By equating the summand of \eqref{v1x3PP} with
the expression of ``$\md{T}_{2n}\pm \md{T}_{2n+1}$" and then canceling
the common factors, we get a polynomial equation
\pq{\label{v1x3YY}\pp{c}{
\md{P}(y)&=
\md{Q}(y)(1-q^{\frac43}y^2)^2(1-q^{\frac32}y^3)(1-q^2y^3)(1-q^{\frac52}y^3)\\
&\pm q^{\frac43}y^3\md{Q}(q^{\frac12}y)(1-q^{\frac12}y)^2(1-q^{\frac56}y)
                     (1-q^{\frac12}y^2)(1-q^{\frac23}y^2),}}
where $q^n$ is replaced by $y$ for brevity.
Now that $\md{P}(y)$ is a polynomial of degree 19 in $y$, we may further take for granted
that $\md{Q}(y)$ is a polynomial of degree 6 so that both sides of the above equation
have the same degree. Writing explicitly
\[\md{Q}(y)=a_0+a_1y+a_2y^2+a_3y^3+a_4y^4+a_5y^5+a_6y^6\]
and then comparing the coefficients of $y^k$ for $0\le k\le 19$ across equation \eqref{v1x3YY},
we get a system of 20 linear equations in seven variables $\{a_0,a_1,a_2,a_3,a_4,a_5,a_6\}$.
The system corresponding to the positive sign `$+$' is not compatible.
Resolving, instead, the system with the minus sign `$-$' of equations
by \emph{Mathematica}, we get the following solution:
\[\Big\{a_0\to1,~a_1\to q^{\frac13},~a_2\to q^{\frac12}+q^{\frac23},~a_3\to0,
~a_4\to-q-q^{\frac56},~a_5\to-q^{\frac76},~a_6\to-q^{\frac32}\Big\}.\]
Therefore, we have determined the polynomial
\pnq{\md{Q}(y)&=1+q^{\frac13} y+q^{\frac12} y^2+q^{\frac23} y^2-q^{\frac56} y^4
-q y^4-q^{\frac76} y^5-q^{\frac32} y^6\\
&=q^{\frac12}y^2(1-q^{\frac13}y^2)(1+q^{\frac16})+(1+q^{\frac13}y)(1-q^{\frac76}y^5).}
Finally, substituting $\md{Q}(y)$ and $\md{T}_n$ into \eqref{v1x3TT}
and simplifying slightly the expression, we arrive at
\pnq{\hyq{c}{q}{q,~q,~q^{\frac56},q^{\frac56}}
            {q^{\frac12},q^{\frac32},q^{\frac13},q^{\frac{4}3}}_\infty
&=\sum_{n=0}^{\infty}(-1)^n
\frac{\ff{q^{\frac12},q^{\frac53};q}{n}
\ff{q^{\frac12},q^{\frac12},q^{\frac13};q^{\frac12}}{n}}
{\ff{q^{\frac13},q^{\frac{4}3};q}{n}(q^{\frac32};q^{\frac12})_{3n}}
q^{\frac{n(9n+7)}{12}}\\
&\times\frac{1-q^{\frac23}}{1-q}
   \bigg\{q^{n+\frac12}
\frac{(1+q^{\frac16})}{(1+q^{\frac{n}2+\frac13})}+
   \frac{1-q^{\frac{5}{2}n+\frac76}}{1-q^{n+\frac13}}
\bigg\},}
which is exactly the identity displayed in Example~\ref{v1x3a}.\qed

The same ``reverse bisection method" can be applied to derive
the formula in Example~\ref{v3x1a}. In fact, by specifying
$\boxed{a=b=c=q^{\frac13},~d=q^{\frac16}}$ in Theorem~\ref{thm=3U},
we can  explicitly write down the formula
\pnq{\hyq{c}{q}{q^{\frac13},q^{\frac13},q^{\frac16},q^{\frac56}}
            {q,~q,~q^{\frac12},~q^{\frac76}}_\infty
&=\sum_{n=0}^{\infty}q^{3n^2+\frac{n}6}
\frac{\ff{q^{\frac13},q^{\frac13},q^{\frac23},q^{\frac16},q^{\frac56},q^{\frac56};q}{n}
\ff{q^{\frac13},q^{\frac16},q^{\frac16};q}{2n}}{\ff{q,q,q^{\frac76};q}{2n}\ff{q,q^{\frac12};q}{3n}}
\mc{W}_n(q^{\frac13},q^{\frac13},q^{\frac13},q^{\frac16}),}
where $\mc{W}_n$ is a rational function of $q^n$ defined by
\pnq{\mc{W}_n(q^{\frac13},q^{\frac13},q^{\frac13},q^{\frac16})
&=\frac{(1-q^{n+\frac23})(1-q^{n+\frac16})(1-q^{2n+\frac16})^2(1-q^{4n+\frac43})}
{(1-q^{2n+1})(1-q^{3n+1})(1-q^{3n+\frac12})}\\
&\times\bigg\{q^{2n+\frac16}
+\frac{(1-q^{2n+1})(1-q^{3n+1})(1-q^{3n+\frac12})(1-q^{4n+\frac16})}
      {(1-q^{n+\frac23})(1-q^{n+\frac16})(1-q^{2n+\frac16})(1-q^{4n+\frac43})}\\
&+\frac{q^{4n+\frac76}(1-q^{n+\frac13})(1-q^{n+\frac56})(1-q^{2n+\frac13})(1-q^{4n+\frac73})}
     {(1-q^{2n+1})(1-q^{3n+2})(1-q^{3n+\frac32})(1-q^{4n+\frac43})}\bigg\}.}
Analogously, it is can be verified without difficulty that
\[\mc{P}(q^n):=\frac{(1-q^{2n+1})^2(1-q^{3n+1})(1-q^{3n+2})(1-q^{3n+\frac12})(1-q^{3n+\frac32})}
{(1-q^{n+\frac16})(1-q^{2n+\frac16})}
\mc{W}_n(q^{\frac13},q^{\frac13},q^{\frac13},q^{\frac16})\]
is a polynomial of degree 19 in $q^n$. Replacing $\mc{W}_n$
by $\mc{P}$, we get the following formula
\pq{\label{v3x1PP}\pp{c}{
\hyq{c}{q}{q^{\frac13},q^{\frac13},q^{\frac56}}
            {q,~q,~q^{\frac32}}_\infty
&=\sum_{n=0}^{\infty}q^{3n^2+\frac{n}6}
\frac{\ff{q^{\frac13},q^{\frac16};q}{2n}
\ff{q^{\frac13};q^{\frac12}}{2n}^2
\ff{q^{\frac16};q^{\frac12}}{2n+1}}
{\ff{q,q;q}{2n+1}\ff{q;q^{\frac12}}{6n+3}}\mc{P}(q^n).}}
According to the factorial structure of the above summand,
we define the hypergeometric term
\[\mc{T}_n:=q^{\frac34n^2+\frac{n}{12}}
\frac{\ff{q^{\frac13},q^{\frac16};q}{n}
\ff{q^{\frac13};q^{\frac12}}{n}^2
\ff{q^{\frac16};q^{\frac12}}{n+1}}
{\ff{q,q;q}{n}\ff{q;q^{\frac12}}{3n}}\mc{Q}(q^{\frac{n}2}),\]
where $\mc{Q}(y)$ is a polynomial to be determined.
Taking into account that
\pnq{\mc{T}_{2n}-\mc{T}_{2n+1}
&=q^{3n^2+\frac{n}6}
\frac{\ff{q^{\frac13},q^{\frac16};q}{2n}
\ff{q^{\frac13};q^{\frac12}}{2n}^2
\ff{q^{\frac16};q^{\frac12}}{2n+1}}
{\ff{q,q;q}{2n+1}\ff{q;q^{\frac12}}{6n+3}}\\
&\times\Big\{\mc{Q}(q^{n})(1-q^{1+2n})^2(1-q^{1+3n})(1-q^{2+3n})(1-q^{\frac32+3n})\\
&-q^{3n+\frac{5}{6}}\mc{Q}(q^{n+\frac12})(1-q^{\frac13+n})^2
(1-q^{\frac23+n})(1-q^{\frac13+2n})(1-q^{\frac16+2n})\Big\},}
we may match the series in \eqref{v3x1PP} with the following ``bisection series"
\pq{\label{v3x1TT}
\hyq{c}{q}{q^{\frac13},q^{\frac13},q^{\frac56}}{q,~q,~q^{\frac32}}_\infty
=\sum_{n=0}^{\infty}(-1)^n\mc{T}_n
=\sum_{n=0}^{\infty}\Big\{\mc{T}_{2n}-\mc{T}_{2n+1}\Big\}.}
By equating the summand of \eqref{v3x1PP} with the expression
of ``$\mc{T}_{2n}-\mc{T}_{2n+1}$" and then canceling
the common factors, we find the polynomial equation below
\pq{\label{v3x1YY}\pp{c}{\mc{P}(y)&=\mc{Q}(y)(1-qy^2)^2(1-qy^3)(1-q^2y^3)(1-q^{\frac32}y^3)\\
&- q^{\frac{5}{6}}y^3\mc{Q}(q^{\frac12}y)(1-q^{\frac13}y)^2(1-q^{\frac23}y)(1-q^{\frac13}y^2)(1-q^{\frac16}y^2).}}
Assuming further that $\mc{Q}(y)$ is of degree 6 and then resolving the linear
system of equations for the coefficients of $\mc{Q}(y)$, we can explicitly determine,
as done for $\md{Q}(y)$ in the precedent example, the following polynomial expression
\pnq{\mc{Q}(y)&=1+q^{\frac16} y+q^{\frac16} y^2+q^{\frac13} y^2-q^{\frac16} y^4
-q^{\frac13} y^4-q^{\frac13} y^5-q^{\frac12} y^6\\
&=q^{\frac16}y^2(1-y^2)(1+q^{\frac16})+(1+q^{\frac16}y)(1-q^{\frac13}y^5).}
By making substitutions of $\mc{Q}(y)$ and $\mc{T}_n$ into \eqref{v3x1TT},
we find, after some simplifications, the reduced formula
as displayed in Example~\ref{v3x1a}:
\pnq{\hyq{c}{q}{q^{\frac43},q^{\frac43},q^{\frac56}}
            {q,~q,~q^{\frac32}}_\infty
&=\sum_{n=0}^{\infty}(-1)^n
\frac{\ff{q^{\frac13},q^{\frac16};q}{n}
\ff{q^{\frac13},q^{\frac13},q^{\frac23};q^{\frac12}}{n}}
{\ff{q,q;q}{n}\ff{q;q^{\frac12}}{3n}}
q^{\frac{n(9n+1)}{12}}\\
&\times\bigg\{q^{n+\frac16}\frac{1-q^{n}}{1-q^{\frac13}}
+\frac{(1+q^{\frac{n}2+\frac16})(1-q^{\frac{5n}2+\frac13})}
{(1+q^{\frac16})(1-q^{\frac13})}\bigg\}.\qqed}

For any given series, it is trivial to write down its bisection series
counterpart. However, reversing this process is quite laborious due
to the complexity to figure out exactly the parameter patterns of the
target series. Notwithstanding the difficulty, we do succeed, with
a help of \emph{Mathematica}, in determining {20} reduced series.
The classical series corresponding to the limits as $q\to1$
with the convergence rate ``$\frac{-1}{27}$" have not been
examined sufficiently, in spite of the fact they are $\pi$-related
series investigated initially by Ramanujan~\cito{ramanujan}
one century ago and explored extensively up to now, for example, in
\cite{kn:adam+w,kn:B+bc,kn:chan-c,kn:chu11bb,kn:chu11mc,kn:chu14mc,kn:chu18e,kn:chu21rj,
	kn:Z-dpv,kn:glaish,kn:gj03em,kn:gj06rj,kn:gj08rj,kn:gj18jdea,kn:zudilin}.

\begin{exam}[$\boxed{a=q^{\frac32},~b=c=d=q}$ in Theorem~\ref{thm=3U}]
\pnq{\hyq{c}{q}{q,~q,~q,~q}
            {q^{\frac12},q^{\frac12},q^{\frac32},q^{\frac32}}_\infty
=\sum_{n=0}^{\infty}(-1)^n
\frac{(q^{\frac12};q^{\frac12})^3_n(q^{\frac12};q)_{n}}
     {(q^{\frac32};q)_{n}(q;q^{\frac12})_{3n}}
q^{\frac{3n(n+1)}4}&\\
\times(1+q^{n+\frac12})
   \bigg\{1+
   \frac{q^{n+\frac12}(1-q^{\frac{n}{2}+\frac12})}
        {(1+q^{n+\frac12})(1-q^{\frac{3}{2}n+1})}
\bigg\}&.}
\end{exam}
When $q\to1$, we recover the following known series 
(cf.~\cite{kn:chu18e,kn:chu14mc}):  
\[\frac{\pi^2}{2}
=\sum_{n=0}^{\infty}\Big(\frac{-1}{27}\Big)^n
\hyp{ccc}{1,\:1}
{\frac43,\frac53\rule[2mm]{0mm}{2mm}}_n
\frac{5+7n}{1+2n}.\]

\begin{exam}[$\boxed{a=q^{\frac12},~b=c=d=q^{\frac13}}$ in Theorem~\ref{thm=3U}]
\pnq{\hyq{c}{q}{q,~q^{\frac13},q^{\frac13},q^{\frac13}}
            {q^{\frac12},q^{\frac16},q^{\frac16},q^{\frac76}}_\infty
=\sum_{n=0}^{\infty}(-1)^n
\frac{(q^{\frac13};q^{\frac12})^3_n(q^{\frac16};q)_{n}}
     {(q^{\frac76};q)_{n}(q^{\frac12};q^{\frac12})_{3n}}
q^{\frac{n(3n+1)}4}&\\
\times(1+q^{n+\frac16})
\bigg\{1+
\frac{q^{n+\frac16}(1-q^{\frac{n}{2}+\frac13})}
     {(1+q^{n+\frac16})(1-q^{\frac{3}{2}n+\frac12})}
\bigg\}&.}
\end{exam}
Its limiting case as $q\to1$ results in the following 
classical series:  
\[\frac{2\pi^2}{\Gam^3(\frac23)}
=\sum_{n=0}^{\infty}\Big(\frac{-1}{27}\Big)^n
\hyp{ccccccccc}{\frac23,\:\frac23,\:\frac16}
{1,\:\frac43,\:\frac76\rule[2mm]{0mm}{2mm}}_n
\big\{8+21n\big\}.\]

\begin{exam}[$\boxed{a=q^{\frac52},~b=q^{\frac23},~c=d=q^{\frac53}}$ in Theorem~\ref{thm=3U}]
\pnq{\hyq{c}{q}{q^2,~q^{\frac23},q^{\frac53},q^{\frac53}}
            {q^{\frac32},q^{\frac56},q^{\frac{5}6},q^{\frac{17}6}}_\infty
=\sum_{n=0}^{\infty}(-1)^n
\frac{(q^{\frac23};q^{\frac12})_n(q^{\frac76};q^{\frac12})^2_n(q^{-\frac16};q)^2_{n}}
     {(q^{\frac{5}6};q)_{n}(q^{\frac{17}6};q)_{n}(q^{\frac32};q^{\frac12})_{3n}}
q^{\frac{n(3n+7)}{4}}&\\
\times\frac{1-q^{\frac16}}{1-q}
   \bigg\{q^{n+\frac56}+
   \frac{(1-q^{\frac{3}{2}n+2})}{(1-q^{n+\frac56})}
   (1+q^{\frac{n}2+\frac16})(1+q^{n+\frac13})
\bigg\}&.}
\end{exam}
Its limiting case as $q\to1$ results in the following 
classical series:  
\[\frac{825\pi^2}{8\Gam^3(\frac13)}
=\sum_{n=0}^{\infty}\Big(\frac{-1}{27}\Big)^n
\hyp{ccccccccc}{\frac{7}3,\:\frac{7}3,-\frac{1}6,-\frac{1}6}
{1,\:\frac53,~\frac{11}6,~\frac{17}6\rule[2mm]{0mm}{2mm}}_n
\big\{53+42 n\big\}.\]

\begin{exam}[$\boxed{a=q^{\frac23},~b=q^{\frac16},~c=d=q}$ in Theorem~\ref{thm=3U}]
\pnq{\hyq{c}{q}{q,~q,~q^{\frac16},~q^{\frac16}}
            {q^{\frac12},q^{\frac32},q^{-\frac13},q^{\frac23}}_\infty
=\sum_{n=0}^{\infty}(-1)^n
\frac{\ff{q^{\frac12},q,q^{-\frac13};q^{\frac12}}{n}
       (q^{\frac12};q)_{n}(q^{\frac43};q)_{n}}
     {(q^{-\frac13};q)_{n}(q^{\frac23};q)_{n}(q^{\frac32};q^{\frac12})_{3n}(1+q^{\frac12})}&\\
\times q^{\frac{n(9n-1)}{12}}
   \bigg\{q^{n+\frac12}+
   \frac{(1-q^{2n+1})(1-q^{\frac{3}2n+\frac16})}
        {(1-q^{\frac{n}2+\frac12})(1-q^{n-\frac13})}
\bigg\}&.}
\end{exam}
Its limiting case as $q\to1$ results in the following 
classical series:  
\[\frac{3\pi\Gamma^2 (\frac{2}{3}) }{\Gamma^2(\frac{1}{6})}
=\sum_{n=1}^{\infty}
\Big(\frac{-1}{27}\Big)^n
\hyp{cccccc}{1,\:\frac{1}2,-\frac{2}3}
{\frac23,\:\frac{2}3,\:\frac{5}3\rule[2mm]{0mm}{2mm}}_n
n\big\{13+21 n\big\}.\]

\begin{exam}[$\boxed{a=q^{\frac13},~b=q^{-\frac16},~c=d=q}$ in Theorem~\ref{thm=3U}]
\pnq{\hyq{c}{q}{q,~q,~q^{-\frac16},q^{-\frac16}}
            {q^{\frac12},q^{\frac32},q^{\frac13},q^{-\frac23}}_\infty
=\sum_{n=0}^{\infty}(-1)^n
\frac{\ff{q^{\frac12},q,q^{-\frac23};q^{\frac12}}{n}
(q^{\frac12};q)_{n}(q^{\frac53};q)_{n}}
     {(q^{-\frac23};q)_{n}(q^{\frac13};q)_{n}(q^{\frac32};q^{\frac12})_{3n}(1+q^{\frac12})}&\\
\times q^{\frac{n(9n-5)}{12}}
   \bigg\{q^{n+\frac12}+
   \frac{(1-q^{2n+1})(1-q^{\frac{3}2n-\frac16})}
    {(1-q^{\frac{n}2+\frac12})(1-q^{n-\frac23})}
\bigg\}&.}
\end{exam}
Its limiting case as $q\to1$ results in the following 
classical series:  
\[\frac{\pi\Gamma^2 (\frac{1}{3}) }{12\Gamma^2(\frac{5}{6})}
=\sum_{n=0}^{\infty}
\Big(\frac{-1}{27}\Big)^n
\hyp{cccccc}{1,\:\frac{1}2,-\frac{4}3}
{\frac13,\:\frac13,\:\frac43\rule[2mm]{0mm}{2mm}}_n
\big\{21 n^2+8 n-3\big\}.\]

\begin{exam}[$\boxed{a=q^{\frac43},~b=c=q,~d=q^{\frac56}}$ in Theorem~\ref{thm=3U}]\label{v1x3a}
\pnq{\hyq{c}{q}{q,~q,~q^{\frac56},q^{\frac56}}
            {q^{\frac12},q^{\frac32},q^{\frac13},q^{\frac{4}3}}_\infty
&=\sum_{n=0}^{\infty}(-1)^n
\frac{(q^{\frac12},q^{\frac53};q)_{n}
\ff{q^{\frac12},q^{\frac12},q^{\frac13};q^{\frac12}}{n}}
     {\ff{q^{\frac13},q^{\frac{4}3};q}{n}(q^{\frac32};q^{\frac12})_{3n}}
q^{\frac{n(9n+7)}{12}}\\
&\times\frac{1-q^{\frac23}}{1-q}
   \bigg\{q^{n+\frac12}
\frac{(1+q^{\frac16})}{(1+q^{\frac{n}2+\frac13})}+
   \frac{1-q^{\frac{5}{2}n+\frac76}}{1-q^{n+\frac13}}
\bigg\}.}
\end{exam}
Its limiting case as $q\to1$ results in the following 
classical series:  
\[\frac{\pi\Gamma^2 (\frac{1}{3}) }{6\Gamma^2(\frac{5}{6})}
=\sum_{n=0}^{\infty}
\Big(\frac{-1}{27}\Big)^n
\hyp{cccccc}{1,\:\frac{1}2,\:\frac{2}3}
{\frac43,\:\frac43,\:\frac43\rule[2mm]{0mm}{2mm}}_n
\big\{3+7 n\big\}.\]

\begin{exam}[$\boxed{a=q^{\frac53},~b=q^{\frac76},~c=d=q}$ in Theorem~\ref{thm=3U}]\label{v1x3b}
\pnq{\hyq{c}{q}{q,~q,~q^{\frac76},~q^{\frac76}}
            {q^{\frac12},q^{\frac32},q^{\frac23},q^{\frac53}}_\infty
=\sum_{n=0}^{\infty}(-1)^n
\frac{\ff{q^{\frac12},q,q^{\frac23};q^{\frac12}}{n}
      (q^{\frac12};q)_{n}(q^{\frac13};q)_{n}}
     {(q^{\frac23};q)_{n}(q^{\frac53};q)_{n}(q^{\frac32};q^{\frac12})_{3n}(1+q^{\frac12})}&\\
\times q^{\frac{n(9n+11)}{12}}
   \bigg\{q^{n+\frac12}+
  \frac{(1-q^{2n+1})(1-q^{\frac{3}2n+\frac76})}
       {(1-q^{\frac{n}2+\frac12})(1-q^{n+\frac23})}
\bigg\}&.}
\end{exam}
Its limiting case as $q\to1$ results in the following 
classical series:  
\[\frac{48\pi\Gamma^2 (\frac{2}{3}) }{\Gamma^2(\frac{1}{6})}
=\sum_{n=0}^{\infty}
\Big(\frac{-1}{27}\Big)^n
\hyp{cccccc}{1,\:\frac{1}2,\:\frac{1}3}
{\frac53,\:\frac{5}3,\:\frac{5}3\rule[2mm]{0mm}{2mm}}_n
\big\{9+28 n+21 n^2\big\}.\]

\begin{exam}[$\boxed{a=b=c=q^{\frac13},~d=q^{\frac16}}$ in Theorem~\ref{thm=3U}]\label{v3x1a}
\pnq{\hyq{c}{q}{q^{\frac43},q^{\frac43},q^{\frac56}}
            {q,~q,~q^{\frac32}}_\infty
&=\sum_{n=0}^{\infty}(-1)^n
\frac{\ff{q^{\frac13},q^{\frac16};q}{n}
\ff{q^{\frac13},q^{\frac13},q^{\frac23};q^{\frac12}}{n}}
{\ff{q,q;q}{n}\ff{q;q^{\frac12}}{3n}}
q^{\frac{n(9n+1)}{12}}\\
&\times\bigg\{q^{n+\frac16}\frac{1-q^{n}}{1-q^{\frac13}}
+\frac{(1+q^{\frac{n}2+\frac16})(1-q^{\frac{5n}2+\frac13})}
{(1+q^{\frac16})(1-q^{\frac13})}\bigg\}.}
\end{exam}
When $q\to1$, we recover the following elegant series 
(cf.~\cite{kn:chu18e,kn:chu21rj,kn:chu14mc}):  
\[\frac{9\sqrt3}{2^{\frac43}\pi}
=\sum_{n=0}^{\infty}
\Big(\frac{-1}{27}\Big)^n
\hyp{cccccc}{\frac13,\:\frac23,\:\frac16}
{1,\:1,\:1\rule[2mm]{0mm}{2mm}}_n
\big\{2 + 21n\big\}.\]

\begin{exam}[$\boxed{a=b=c=q^{\frac23},~d=q^{\frac56}}$ in Theorem~\ref{thm=3U}]\label{v3x1b}
\pnq{\hyq{c}{q}{q^{\frac{2}3},q^{\frac{5}3},q^{\frac16}}
            {q,~q,~q^{\frac12}}_\infty
=\sum_{n=0}^{\infty}(-1)^n
\frac{\ff{q^{\frac16},q^{\frac16},q^{\frac13};q^{\frac12}}{n}(q^{\frac53};q)_{n}(q^{\frac56};q)_{n}}
     {(q;q)^2_{n}(q^{\frac12};q^{\frac12})_{3n}
(1+q^{\frac{n}2+\frac13}+q^{n+\frac23})}&\\
\times q^{\frac{9n^2+5n+4}{12}}
\bigg\{\frac{1-q^{n-\frac16}}{1-q^{\frac{n}2+\frac13}}
+\frac{q^{-\frac{n}2-\frac13}(1-q^{2n+\frac23})^2}{(1-q^{n+\frac23})(1-q^{\frac{3n}2+\frac12})}
\bigg\}&.}
\end{exam}
When $q\to1$, we recover another elegant series 
(cf.~\cite{kn:chu18e,kn:chu21rj,kn:chu14mc}):  
\[\frac{27\sqrt3}{2^{\frac53}\pi}
=\sum_{n=0}^{\infty}
\Big(\frac{-1}{27}\Big)^n
\hyp{cccccc}{\frac13,\:\frac23,\:\frac56}
{1,\:1,\:1\rule[2mm]{0mm}{2mm}}_n
\big\{5 + 42n\big\}.\]

\begin{exam}[$\boxed{a=b=q^{\frac53},~c=q^{\frac23},~d=q^{\frac56}}$ in Theorem~\ref{thm=3U}]\label{v3x1c}
\pnq{\hyq{c}{q}{q^{\frac53},q^{\frac53},q^{\frac76},q^{\frac76}}
            {q^2,~q^2,q^{\frac32},q^{\frac16}}_\infty
=\sum_{n=0}^{\infty}(-1)^n
\frac{\ff{q^{\frac23},q^{\frac56},q^{\frac76};q^{\frac12}}{n}(q^{-\frac16};q)_{n}(q^{\frac23};q)_{n}}
     {(q;q)_{n}(q^2;q)_{n}(q^{\frac32};q^{\frac12})_{3n}(1+q^{\frac13})}
q^{\frac{n(9n+11)}{12}}&\\
\times
   \bigg\{(1+q^{\frac{n}{2}+\frac16})(1+q^{\frac{n}2+\frac13})(1+q^{\frac{n}2+\frac12})(1+q^{n+\frac13})
-q^{\frac{n+1}2}(1+q^{\frac{3}{2}n+\frac13})
\bigg\}&.}
\end{exam}
Its limiting case as $q\to1$ results in the following 
classical series:  
\[\frac{81\sqrt3}{28\cdot2^{2/3}\pi}
=\sum_{n=0}^{\infty}
\hyp{ccc}{\frac23,\frac73,-\frac16}
{1,~1,~2\rule[-2mm]{0mm}{2mm}}_n
\Big(\frac{-1}{27}\Big)^n.\]

\begin{exam}[$\boxed{a=q^{\frac43},~b=q^{\frac13},~c=q^{\frac43},~d=q^{\frac76}}$ in Theorem~\ref{thm=3U}]
\pnq{\hyq{c}{q}{q^{\frac13},q^{\frac43},q^{\frac76},q^{\frac{11}6}}
               {q,~q^2,~q^{\frac32},~q^{\frac16}}_\infty
=\sum_{n=0}^{\infty}(-1)^n
\frac{\ff{q^{\frac23},q^{\frac43},q^{-\frac16};q^{\frac12}}{n}
      (q^{\frac13};q)_{n}(q^{\frac76};q)_{n}}
     {(q;q)^2_{n}(q^{\frac32};q^{\frac12})_{3n}}&\\
\times q^{\frac{n(9n+7)}{12}}
   \bigg\{q^{n+\frac16}+
  \frac{(1-q^{\frac{3}2n+\frac56})1-q^{2n+\frac43})}
       {(1-q^{\frac{n}2+\frac56})(1-q^{n+1})}\bigg\}&.}
\end{exam}
Its limiting case as $q\to1$ results in the following 
classical series:  
\[\frac{81\sqrt3}{2^{\frac13}\pi}
=\sum_{n=0}^{\infty}
\Big(\frac{-1}{27}\Big)^n
\hyp{cccccc}{\frac13,-\frac13,\:\frac76}
{1,~1,~2\rule[2mm]{0mm}{2mm}}_n
\big\{35+90 n+63 n^2\big\}.\]

\subsection{} \
The $\boxed{\del=1}$ case of Theorem~\ref{thm=V} under the partition
pattern \eqref{p3pattern} results in the formula below.
\begin{thm}[Nonterminating series identity]\label{thm=3V}
\pnq{\hyq{c}{q}{b,~c,~qd,~qa^2/bcd}{qa/b,qa/c,a/d,bcd/a}_\infty
=1+\sum_{n=1}^{\infty}W_n(a,b,c,d)
\Big(\frac{q^{3n-2}a^3}{bcd}\Big)^n&\\
\times\frac{\ff{b,c,d,qa/bc,qa/bd,a/cd;q}{n}\ff{bc/a,bd/a,cd/a;q}{2n}}
{\ff{a/b,a/c,a/d;q}{2n}\ff{q,bcd/a;q}{3n}}&}
\text{where the weight function is given by}
\pnq{W_n(a,b,c,d)
&=\tfrac{(1-q^{3n})(1-b/a)(1-c/a)(1-q^{3n-1}bcd/a) (1-q^{4n-2}b)}
      {(1-cd/a)(1-1/d)(1-q^{n-1}b)(1-q^{2n-1}bc/a)(1-q^{2n-1}bd/a)}\\
\times\bigg\{1&-
\tfrac{(1-q^{n-1}b)(1-q^{2n-1}bc/a) (1-q^{2n-1}bd/a)(1-q^{n}a/(cd))(1-q^{4n}d)}
{(1-q^{3n})(1-q^{-2n}b/a) (1-q^{2n}a/c)(1-q^{3n-1}bcd/a)(1-q^{4n-2}b)}\\
&-\tfrac{(1-q^{3n-1})(1-q^{1-2n}c/a)(1-q^{2n-1}a/d)(1-q^{3n-2}bcd/a)(1-q^{4n-3}c)}
{(1-q^{n-1}c)(1-q^{2n-2}bc/a)(1-q^{n}a/(bd))(1-q^{2n-1}cd/a)(1-q^{4n-2}b)}\bigg\}.}
\end{thm}

Following the same ``reverse bisection approach" illustrated
for Examples~\ref{v1x3a} and ~\ref{v3x1a}, we can derive
from Theorem~\ref{thm=3V} further {9} identities.
\begin{exam}[$\boxed{a=q^{\frac32},~b=c=d=q}$ in Theorem~\ref{thm=3V}]
\pnq{\hyq{c}{q}{q~,q~,q,~q^2}
            {q^{\frac12},q^{\frac32},q^{\frac32},q^{\frac32}}_\infty
=1+\sum_{n=0}^{\infty}(-1)^n
\frac{\ff{q^{\frac12},q^{\frac12},q^{\frac32};q^{\frac12}}{n}(q^{\frac12};q)_{n}}
     {(q^{\frac52};q)_{n}(q^{2};q^{\frac12})_{3n}}
q^{\frac{(n+2)(3n+1)}4}&\\
\times\frac{(1-q^{\frac12})(1-q^{2n+2})}
      {(1+q^{\frac12})(1-q^{\frac32})^2}
\bigg\{1+\frac{q^{n+\frac32}(1-q^{\frac{n}2+\frac12})(1-q^{n+\frac12})}
{(1-q^{\frac{3}2n+2})(1-q^{2n+2})}\bigg\}&.}
\end{exam}
Its limiting case as $q\to1$ results in the following 
classical series:  
\[9\pi^2
=89+\sum_{n=1}^{\infty}\Big(\frac{-1}{27}\Big)^n
\hyp{ccccccccc}{\frac12,\:1,\:3}
{\frac52,\:\frac53,\:\frac73\rule[2mm]{0mm}{2mm}}_n
\big\{17+14n\big\}.\]

\begin{exam}[$\boxed{a=q^{\frac12},~b=c=q^{\frac23},~d=q^{-\frac13}}$ in Theorem~\ref{thm=3V}]
\pnq{\hyq{c}{q}{q,~q^{\frac23},q^{\frac23},q^{\frac23}}
            {q^{\frac12},q^{\frac56},q^{\frac56},q^{\frac56}}_\infty
=&\sum_{n=0}^{\infty}(-1)^n
\frac{(q^{\frac16};q^{\frac12})^2_n(q^{\frac23};q^{\frac12})(q^{-\frac16};q)_{n}}
     {(q^{\frac56};q)_{n}(q^{\frac12};q^{\frac12})_{3n}}\\
\times&q^{\frac{3n(n+1)}{4}}
\bigg\{1+\frac{(1-q^{\frac{3}2n})(1+q^{\frac16-n})}
{1-q^{\frac{n}2+\frac16}}\bigg\}.}
\end{exam}
Its limiting case as $q\to1$ results in the following   
classical series:  
\[\frac{2\pi^2}{\Gam^3(\frac13)}
=\sum_{n=0}^{\infty}\Big(\frac{-1}{27}\Big)^n
\hyp{ccc}{\frac{1}3,\:\frac{1}3}
{1,~\frac23\rule[2mm]{0mm}{2mm}}_n
\frac{1+21 n}{1-6n}.\]

\begin{exam}[$\boxed{a=q^{\frac32},~b=q^{\frac43},~c=d=q^{\frac13}}$ in Theorem~\ref{thm=3V}]
\pnq{\hyq{c}{q}{q^2,~q^{\frac13},q^{\frac43},q^{\frac43}}
            {q^{\frac12},q^{\frac16},q^{\frac76},q^{\frac{13}6}}_\infty
&=\sum_{n=0}^{\infty}(-1)^n
\frac{\ff{q^{\frac13},q^{\frac43},q^{\frac56};q^{\frac12}}{n}(q^{-\frac56};q)^2_{n}}
     {(q^{\frac76};q)_{n}(q^{\frac{13}6};q)_{n}(q^{\frac12};q^{\frac12})_{3n}}
q^{\frac{3n^2+3n+2}{4}}\\
&\times\frac{1+q^{\frac{n}2-\frac23}}
{1-q^{\frac56}}
\bigg\{q^{n+\frac{1}6}
+\frac{(1-q^{2n+\frac56})(1+q^{n-\frac13}+q^{\frac{n}2-\frac16})}
{(1-q^{\frac{n}2+\frac56})(1+q^{\frac{n}2-\frac23})}\bigg\}.}
\end{exam}
Its limiting case as $q\to1$ results in the following 
classical series:  
\[\frac{175\pi^2}{36\Gam^3(\frac23)}
=\sum_{n=0}^{\infty}\Big(\frac{-1}{27}\Big)^n
\hyp{ccccccccc}{\frac53,\:\frac53,-\frac56,-\frac56}
{1,~\frac13,~\frac{7}6,~\frac{13}6\rule[2mm]{0mm}{2mm}}_n
\big\{25+42 n\big\}.\]

\begin{exam}[$\boxed{a=q^{\frac43},~b=q^{\frac56},~c=d=q}$ in Theorem~\ref{thm=3V}]
\pnq{\hyq{c}{q}{q,~q^2,~q^{\frac56},q^{\frac{5}6}}
            {q^{\frac32},q^{\frac32},q^{\frac13},q^{\frac43}}_\infty
=1+\sum_{n=0}^{\infty}(-1)^{n}
\frac{\ff{q^{\frac12},q^{\frac32},q^{\frac56};q^{\frac12}}{n}
       \ff{q^{\frac12},q^{\frac23};q}{n}}
     {(q^{\frac43};q)^2_{n}(q^{\frac32};q^{\frac12})_{3n}}
q^{\frac{9n^2+19n+4}{12}}&\\
\times\frac{(1-q^{\frac12})(1+q^{\frac{n}2+\frac13})}
{(1+q^{\frac12})(1-q^{2+\frac{3n}2})}
\bigg\{q^{\frac{n+1}2}+\frac{(1-q^{2n+\frac{11}6})(1+q^{n+\frac32})}
{(1-q^{n+\frac43})(1+q^{n+1}+q^{\frac{n+1}2})}
\bigg\}&.}
\end{exam}
Its limiting case as $q\to1$ results in the following 
classical series:  
\[\frac{4\pi\Gamma^2 (\frac{1}{3}) }{\Gamma^2(\frac{5}{6})}
=71+\sum_{n=1}^{\infty}
\Big(\frac{-1}{27}\Big)^n
\hyp{cccccc}{3,\:\frac{1}2,\:\frac{2}3}
{\frac43,\:\frac{7}3,\:\frac{7}3\rule[2mm]{0mm}{2mm}}_n
\big\{23+21n\big\}.\]

\begin{exam}[$\boxed{a=q^{\frac23},~b=q^{\frac16},~c=d=q}$ in Theorem~\ref{thm=3V}]
\pnq{\hyq{c}{q}{q,~q^2,~q^{\frac16},q^{\frac{1}6}}
            {q^{\frac32},q^{\frac32},q^{-\frac13},q^{\frac23}}_\infty
=1-\sum_{n=0}^{\infty}(-1)^{n}
\frac{\ff{q^{\frac12},q^{\frac32},q^{-\frac13};q^{\frac12}}{n}
       \ff{q^{\frac12},q^{\frac{7}3};q}{n}}
     {(q^{\frac23};q)^2_n(q^{\frac32};q^{\frac12})_{3n}}
q^{\frac{n(9n+11)}{12}}&\\
\times\frac{(1-q^{\frac12})(1+q^{\frac23})(1-q^{\frac23})}
{(1+q^{\frac12})(1-q^{\frac13})(1-q^{2+\frac{3n}2})}
\bigg\{q^{\frac{n+1}2}+\frac{(1-q^{2n+\frac76})(1+q^{n+\frac32})}
{(1-q^{n+\frac23})(1+q^{n+1}+q^{\frac{n+1}2})}\bigg\}&.}
\end{exam}
Its limiting case as $q\to1$ results in the following 
classical series:  
\[\frac{9\pi\Gamma^2 (\frac{2}{3}) }{\Gamma^2(\frac{1}{6})}
=1+\sum_{n=1}^{\infty}
\Big(\frac{-1}{27}\Big)^n
\hyp{cccccc}{3,\:\frac{1}2,-\frac{2}3}
{\frac23,\:\frac{5}3,\:\frac{5}3\rule[2mm]{0mm}{2mm}}_n
\big\{13+21n\big\}.\]

\begin{exam}[$\boxed{a=b=c=q^{\frac23},~d=q^{-\frac16}}$ in Theorem~\ref{thm=3V}]\label{v3x1d}
\pnq{\hyq{c}{q}{q^{\frac23},q^{\frac23},q^{\frac{7}6}}
            {q,~q,~q^{\frac12}}_\infty
=&\sum_{n=0}^{\infty}(-1)^{n}
\frac{(q^{\frac16};q^{\frac12})^2_n(q^{\frac56};q^{\frac12})_n
       \ff{q^{-\frac13},q^{\frac56};q}{n}}
     {(q;q)^2_{n}(q^{\frac12};q^{\frac12})_{3n}}\\
\times&q^{\frac{n(9n+11)}{12}}
\bigg\{1+q^{\frac16-n}
\frac{(1-q^{\frac{3}2n})(1-q^{2n-\frac16})}
{(1-q^{\frac{n}2+\frac13})(1-q^{n-\frac16})}\bigg\}.}
\end{exam}
Its limiting case as $q\to1$ results in the following 
classical series:  
\[\frac{3\sqrt3}{\pi\sqrt[3]4}
=\sum_{n=0}^{\infty}
\Big(\frac{-1}{27}\Big)^n
\hyp{cccccc}{\frac13,-\frac13,-\frac16}
{1,~1,~1\rule[2mm]{0mm}{2mm}}_n
\big\{1-63 n^2\big\}.\]

\begin{exam}[$\boxed{a=q^{\frac13},~b=q^{\frac16},~c=d=q^{\frac13}}$ in Theorem~\ref{thm=3V}]
\pnq{\hyq{c}{q}{q^{\frac43},q^{\frac43},q^{\frac{11}6}}
            {q,~q,~q^{\frac52}}_\infty
&=\sum_{n=0}^{\infty}(-1)^n
\frac{\ff{q^{\frac13},q^{\frac43},q^{\frac16};q^{\frac12}}{n}
       \ff{q^{\frac43},q^{\frac16};q}{n}}
     {(q;q)_{n}(q;q)_{n}(q^2;q^{\frac12})_{3n}}
q^{\frac{n(9n+13)}{12}}\\
&\times\frac{1-q^{\frac{3n}2+1}}{1-q}
\bigg\{q^{\frac{n}2+\frac13}
+\frac{(1-q^{2n+\frac76})(1+q^{n+\frac76})}
{(1-q^{n+1})(1+q^{\frac{n}2+\frac13}+q^{n+\frac23})}\bigg\}.}
\end{exam}
Its limiting case as $q\to1$ results in the following 
classical series:  
\[\frac{729\sqrt3}{20\sqrt[3]{2}\pi}
=\sum_{n=0}^{\infty}
\Big(\frac{-1}{27}\Big)^n
\hyp{cccccc}{\frac13,\:\frac83,\:\frac16}
{1,\:2,\:2\rule[2mm]{0mm}{2mm}}_n
\big\{16 + 21n\big\}.\]

\begin{exam}[$\boxed{a=q^{\frac23},~b=q^{\frac56},~c=d=q^{\frac23}}$ in Theorem~\ref{thm=3V}]
\pnq{\hyq{c}{q}{q^{\frac53},q^{\frac53},q^{\frac{11}6},q^{\frac{13}6}}
            {q,~q^3,~q^{\frac52},q^{\frac56}}_\infty
=&\sum_{n=0}^{\infty}(-1)^n
\frac{\ff{q^{\frac16},q^{\frac13},q^{\frac53};q^{\frac12}}{n}
       \ff{q^{\frac53},q^{\frac{11}6};q}{n}}
     {(q;q)_{n}(q;q)_{n}(q^{\frac52};q^{\frac12})_{3n}}\\
\times&q^{\frac{(n+2)(9n+5)}{12}}
\bigg\{1+\frac{(1-q^{2n+\frac83})(1-q^{\frac{3}2n+\frac76})}
{q^{n+\frac56}(1-q^{\frac{n}2+\frac76})(1-q^{n+1})}\bigg\}.}
\end{exam}
Its limiting case as $q\to1$ results in the following 
classical series:  
\[\frac{2187 \sqrt3}{10\sqrt[3]{4}\pi}
=\sum_{n=0}^{\infty}
\Big(\frac{-1}{27}\Big)^n
\hyp{cccccc}{\frac13,\:\frac{2}3,\:\frac{11}6}
{1,~2,~2\rule[2mm]{0mm}{2mm}}_n
\big\{77+144 n+63 n^2\big\}.\]

\begin{exam}[$\boxed{a=q^{\frac43},~b=q^{\frac76},~c=q^{\frac13},~d=q^{\frac43}}$ in Theorem~\ref{thm=3V}]
\pnq{\hyq{c}{q}{q^{\frac73},q^{\frac73},q^{\frac{11}6}}
            {q^2,~q^2,~q^{\frac52}}_\infty
=&\sum_{n=0}^{\infty}(-1)^n
\frac{\ff{q^{\frac23},q^{-\frac16},q^{\frac{11}6};q^{\frac12}}{n}
       \ff{q^{\frac73},q^{\frac16};q}{n}}
     {(q;q)_{n}(q^2;q)_{n}(q^{\frac52};q^{\frac12})_{3n}}
q^{\frac{n(9n+19)}{12}}\\
\times&\bigg\{\frac{q^{n+\frac43}(1-q^{n})}
      {(1-q^2)(1+q^{\frac{n}2+\frac23})}
+\frac{(1-q^{\frac{3}2n+\frac{11}6})(1-q^{2n+\frac{13}6})}
      {(1-q^2)(1-q^{n+2})}\bigg\}.}
\end{exam}
Its limiting case as $q\to1$ results in the following 
classical series:  
\[\frac{6561 \sqrt3}{40\sqrt[3]{2}\pi}
=\sum_{n=0}^{\infty}
\Big(\frac{-1}{27}\Big)^n
\hyp{cccccc}{\frac{4}3,\frac{11}3,\frac{1}6}
{1,\:~2,\:~3\rule[2mm]{0mm}{2mm}}_n
\frac{143+285 n+126 n^2}{(1-3n)(2+3n)}.\]

From the above examples, one may get an impression that all the series
obtained through triplicate inversions might be reducible, via the
``reverse bisection method", to simpler ones. However, this is not
true  in general. Five counterexamples below demonstrate this fact.
\begin{exam}[$\boxed{a=b=c=q^{\frac12},~d=q^{\frac14}}$ in Theorem~\ref{thm=3U}]
\pnq{\hyq{c}{q}{q^{\frac12},q^{\frac52}}
            {q,~q^2}_\infty
=&\sum_{n=0}^{\infty}q^{3n^2+\frac{n}{4}}
\frac{(q^{\frac12};q)^2_n(q^{\frac32};q)_n(q^{\frac34};q)^2_{n}(q^{\frac14};q)_{n+1}
      (q^{\frac12};q)_{2n}(q^{\frac14};q)_{2n+1}}
{(q;q)_{2n}(q^2;q)_{2n}(q;q)_{3n+1}(q^{\frac34};q)_{3n+1}}\\
\times\frac{1-q^{4n+\frac32}}{1-q^{\frac32}}
&\bigg\{q^{2n+\frac14}+
\frac{(1-q^{2n+1})(1-q^{3n+1})(1-q^{3n+\frac34})(1-q^{4n+\frac14})}
      {(1-q^{n+\frac14})(1-q^{n+\frac12})(1-q^{2n+\frac14})(1-q^{4n+\frac32})}\\
&+\frac{q^{4n+\frac54}(1-q^{n+\frac12})(1-q^{n+\frac34})(1-q^{2n+\frac12})(1-q^{4n+\frac52})}
     {(1-q^{2n+1})(1-q^{3n+\frac74})(1-q^{3n+2})(1-q^{4n+\frac32})}\bigg\}.}
\end{exam}
Its limiting case as $q\to1$ results in the following 
classical series:  
\pnq{\frac{16}{\pi }&=\sum_{n=0}^{\infty}\frac{1}{729^n}
\hyp{cccccc}{\frac12,\frac14,\frac34,\frac34,\frac34,\frac58,\frac98,\frac{11}8}
{1,1,1,\frac23,\frac43,\frac38,\frac{7}{12},\frac{11}{12} \rule[2mm]{0mm}{2mm}}_n\\
&\times\bigg\{1+\frac{12 (3 n+1) (16 n+1)}{(8 n+1) (8 n+3)}
+\frac{(4 n+1) (4 n+3) (8 n+5)}{4 (3 n+2) (8 n+3) (12 n+7)}\bigg\}.}

\begin{exam}[$\boxed{a=q^{\frac32},~b=q^{\frac56},~c=d=q}$ in Theorem~\ref{thm=3V}]
\pnq{\hyq{c}{q}{q,~q^2,q^{\frac56},q^{\frac76}}
            {q^{\frac12},q^{\frac32},q^{\frac43},q^{\frac53}}_\infty
=1+\sum_{n=1}^{\infty}q^{\frac{9n^2+11n}{3}}
\frac{(q^{-\frac12};q)_n(q;q)^2_n(q^{\frac23};q)^2_{n}(q^{-\frac16};q)_{n}(q^{-\frac23};q)^2_{2n}}
     {(q^{\frac12};q)_{2n}(q^{\frac23};q)_{2n}(q;q)_{3n}(q^{\frac13};q)_{3n}}&\\
\times\frac{(1-q^{3n})(1-q^{\frac76-4n})}{(1-q)(1+q^{\frac13})(1-q^{-\frac16})}
\Bigg\{1-\frac{(1-q^{\frac12-2n})(1-q^{3n-1})(1-q^{3n-\frac23})(1-q^{4n-2})}
{(1-q^{n})(1-q^{n-\frac13})(1-q^{2n-\frac53})(1-q^{4n-\frac76})}&\\
+\frac{q^{2n+\frac23}(1-q^{n-\frac12})(1-q^{n-\frac16})(1-q^{2n-\frac23})^2(1-q^{4n+1})}
{(1-q^{3n})(1-q^{2n+\frac12})(1-q^{2n+\frac23})(1-q^{3n+\frac13})(1-q^{4n-\frac76})}\Bigg\}&.}
\end{exam}
Its limiting case as $q\to1$ results in the following 
classical series:  
\pnq{\frac{2\pi}{3\sqrt{3}}
&=1+\sum_{n=1}^{\infty}\frac{1}{729^n}
\hyp{cccccc}{1,\frac{1}{2},\frac{2}{3},\frac{2}{3},\frac{2}{3},~\frac{1}{6},\frac{1}{6}}
{\frac13,\frac43,\frac14,\frac34,\frac{4}{9},\frac{7}{9},\frac{10}{9} \rule[2mm]{0mm}{2mm}}_n\\
&\times\bigg\{1+\frac{3 n (3 n+1) (9 n+1) (24 n-7)}{2 (2 n-1) (3 n-1)^2 (6n-1)}
+\frac{27 (3 n+1) (4 n-1) (9 n-2) (9 n+1)}{(3 n-1)^2 (6 n-5) (6 n-1)}\bigg\}.}

\begin{exam}[$\boxed{a=b=c=d=q^{\frac12}}$ in Theorem~\ref{thm=3U}]\label{w1+1+1a}
\pnq{\hyq{c}{q}{q^{\frac12},q^{\frac32},q^{\frac32},q^{\frac32}}
            {q,~q,~q^2,~q^2}_\infty
&=\sum_{n=0}^{\infty}q^{3n^2}
\frac{(q^{\frac12};q)_n^6(q^{\frac32};q)^3_{2n}(1-q^{n+\frac12})^2(1-q^{4n+\frac32})}
{(q;q)^3_{2n}(q^2;q)^2_{3n}(1-q^{2n+1})(1-q^{2n+\frac12})}\\
&\times\bigg\{q^{2n}
\pp{t}{&+\frac{(1+q^{n+\frac12})(1-q^{2n})(1-q^{3n+1})^2(1-q^{4n+\frac12})}
      {(1-q^{n+\frac12})(1-q^{2n+\frac12})^2(1-q^{4n+\frac32})}\\
&+\frac{q^{4n+1}(1-q^{2n+\frac12})(1-q^{2n+\frac32})(1-q^{4n+\frac52})}
     {(1+q^{n+\frac12})^2(1-q^{3n+2})^2(1-q^{4n+\frac32})}\bigg\}.}}
\end{exam}
Its limiting case as $q\to1$ results in the following 
classical series:  
\pnq{\frac{32}{\pi^2}&=\sum_{n=0}^{\infty}
\frac{(\frac12)^6_n(\frac12)^3_{2n}}{(2n)!^3(3n)!^2}
\times\frac{(1+2n)(1+4n)^2(3+8n)}{(1+3n)^2}\\
&\times\Bigg\{1+\frac{32 n (1+3 n)^2 (1+8 n)}{(1+2 n) (1+4 n)^2 (3+8 n)}
+\frac{(1+4 n) (3+4 n) (5+8 n)}{16 (2+3 n)^2 (3+8 n)}\Bigg\}.}

\begin{exam}[$\boxed{a=b=c=d=q^{\frac12}}$ in Theorem~\ref{thm=3V}]\label{w1+1+1b}
\pnq{\hyq{c}{q}{q^{\frac12},q^{\frac52},q^{\frac52},q^{\frac52}}
            {q,~q,~q^3,~q^3}_\infty
&=\sum_{n=0}^{\infty}q^{(n+1)(3n+1)}
\frac{(q^{\frac12};q)^6_n(q^{\frac52};q)^3_{2n}(1-q^{n+\frac12})^4(1-q^{4n+\frac52})}
{(q^2;q)^3_{2n}(q^3;q)^2_{3n}(1-q)^3(1-q^{2n+\frac32})^2}\\
&\times\bigg\{1
\pp{t}{&+\frac{q^{-1-2n}(1-q^{2n+1})^2(1-q^{3n+2})^2(1-q^{4n+\frac32})}
      {(1-q^{n+\frac12})^2(1-q^{2n+\frac12})(1-q^{2n+\frac32})(1-q^{4n+\frac52})}\\
&+\frac{q^{2+2n}(1-q^{n+\frac12})^2(1-q^{2n+\frac32})^2(1-q^{4n+\frac92})}
     {(1-q^{2n+2})^2(1-q^{3n+3})^2(1-q^{4n+\frac52})}\bigg\}.}}
\end{exam}
Its limiting case as $q\to1$ results in the following 
classical series:  
\pnq{\frac{512}{\pi^2}&=\sum_{n=0}^{\infty}
\frac{(\frac12)^6_n(\frac12)^3_{2n}}{(2n)!^3(3n)!^2}
\times\frac{(1+2n)(1+4n)^3(3+4n)(5+8n)}{(1+3n)^2(2+3n)^2}\\
&\times\Bigg\{1+\frac{16 (2+3 n)^2 (3+8 n)}{(1+4 n) (3+4 n)(5+8 n)}
+\frac{(1+2n)^2 (3+4 n)^2 (9+8 n)}{576 (1+ n)^4 (5+8 n)}\Bigg\}.}


These last two $q$-series have the same structure of factorial quotients, and differ
to each other only by extra factors of rational functions in $q^n$. When $q\to1$, their
limiting series with the convergence rate ``$\frac{1}{729}$" are also significant, since
they have not been recorded in the literature, even though there have been extensive
research activities around infinite series for $1/\pi^2$, mainly made
by Guillera~\cite{kn:gj03em,kn:gj06rj,kn:gj08rj}.

\subsection{} \
Theorems~\ref{thm=U} and~\ref{thm=V} are so general that can be specified
concretely to produce numerous summation formulae. For instance, it is well
known that $\boxed{n=\pavv{\tfrac{1+n}3}+\pavv{\tfrac{1+2n}3}}$ holds
for all $n\in\mb{N}_0$. Then it is not difficult to check that Jackson's
formula \eqref{jackson} is equivalent to the following one
\[\Ome_n\big(a;q^{\pav{\tfrac{1+n}3}}b,c,q^{\pav{\tfrac{1+2n}3}}d\big)
{=}\hyq{c}{q}{qa,~bd/a}{qa/c,bcd/a}_n
\hyq{c}{q}{qa/cd,bc/a}{qa/d,b/a}_{\pavv{\frac{1+n}3}}
\hyq{c}{q}{qa/bc,cd/a}{qa/b,d/a}_{\pavv{\frac{1+2n}3}}\]
with its parameters subject to $\boxed{qa^2=bcde}$.
The corresponding partition pattern is determined by
\[\boxed{\Lam=3:\mult{cccc}
{\veps_b=1,&\veps_c=0,&\veps_d=1,&\veps_e=0;\\
\lam_b=1,&\lam_c=0,&\lam_d=2,&\lam_e=0.}}\]
For Theorems~\ref{thm=U} and \ref{thm=V}, there exist six different
summation formulae in accordance with $\del=0,1,2$. We limit to highlight
only the $\boxed{\del=0}$ case of Theorem~\ref{thm=U} together with two
infinite series identities as examples. When $q\to1$, the corresponding
limiting series have the convergence rate ``$\frac{16}{729}$".
\begin{thm}[Nonterminating series identity]\label{thm=p23U}
\pnq{\hyq{c}{q}{b,~c,~qd,~qa^2/bcd}{qa/b,qa/c,a/d,bcd/a}_\infty
=\sum_{n=0}^{\infty}W_n(a,b,c,d)
\Big(\frac{q^{2n}a^3}{bd^2}\Big)^n&\\
\times\frac{\ff{qa/cd,b,bc/a;q}{n}\ff{qa/bc,qd,cd/a;q}{2n}\ff{bd/a;q}{3n}}
{\ff{qa/d;q}{n}\ff{qa/b;q}{2n}\ff{q,qa/c,bcd/a;q}{3n}},&}
where the weight function is given by
\pnq{W_n&(a,b,c,d)
=\frac{(1-q^{2n+1}a/bc)(1-q^{2n}cd/a)(1-q^{3n}bd/a)(1-q^{4n+1}b)}
      {(1-q^{3n+1})(1-q^{3n+1}a/c)(1-a/d)(1-q^{3n}bcd/a)}\\
\times~&\bigg\{q^n\frac{a}{d}+
\frac{(1-q^na/d)(1-q^{3n+1})(1-q^{3n+1}a/c)(1-q^{3n}bcd/a)(1-q^{5n}d)}
     {(1-q^{2n+1}a/bc)(1-q^{2n}cd/a)(1-q^{2n}d)(1-q^{3n}bd/a)(1-q^{4n+1}b)}\\
+~&\frac{q^{3 n+1}a^2 (1-q^nb)(1-q^nbc/a) (1-q^{n+1}a/cd) (1-q^{3n+1}bd/a)(1-q^{5 n+3}d)}
      {bd (1-q^{2 n+1}a/b)(1-q^{3n+2})(1-q^{3n+2}a/c)(1-q^{3n+1}bcd/a)(1-q^{4n+1}b)}\bigg\}.}
\end{thm}

\begin{exam}[$\boxed{a=b=c=d=q^{\frac12}}$ in Theorem~\ref{thm=p23U}]\label{w1+2d}
\pnq{\hyq{c}{q}{q^{\frac12},q^{\frac12},q^{\frac12},q^{\frac52}}{q,~q,~q,~q}_\infty
&=\sum_{n=0}^{\infty}q^{2n^2}
\frac{(q^{\frac12};q)^3_n(q^{\frac12};q)^3_{2n+1}(q^{\frac32};q)_{3n}(1-q^{4n+\frac32})}
     {(q;q)_n(q;q)_{2n}(q;q)^3_{3n+1}(1-q^{3/2})}\\
&\times
   \bigg\{q^n+
\frac{(1-q^n) (1-q^{3n+1})^3 (1-q^{5 n+\frac{1}{2}})}
     {(1-q^{3 n+\frac{1}{2}})(1-q^{2n+\frac{1}{2}})^3(1-q^{4n+\frac{3}{2}})}\\
&\qquad+q^{3n+1}
\frac{(1-q^{n+\frac{1}{2}})^3 (1-q^{3 n+\frac{3}{2}})(1-q^{5 n+\frac{7}{2}})}
     {(1-q^{3 n+2})^3(1-q^{2 n+1})(1-q^{4n+\frac{3}{2}})}\bigg\}.}
\end{exam}
This provides a $q$-analogue for the following series:
\pnq{\frac{32}{\pi^2}
&=\sum_{n=0}^{\infty}(3+8n)
\dfrac{(\frac12)^3_n(\frac32)^3_{2n}(\frac32)_{3n}}
     {n!(2n)!(3n+1)!^3}\\
&\times\bigg\{1+\frac{16 n (1+3 n)^3 (1+10 n)}{(1+4 n)^3 (1+6 n) (3+8 n)}
+\frac{3 (1+2 n)^3 (7+10 n)}{16 (2+3 n)^3 (3+8 n)}\Big\}.}

\begin{exam}[$\boxed{a=q^{\frac32},~b=c=d=q}$ in Theorem~\ref{thm=p23U}]\label{w1+2e}
\pnq{\hyq{c}{q}{q,~q,~q,~q}{q^{\frac12},q^{\frac12},q^{\frac12},q^{\frac32}}_\infty
&=\sum_{n=0}^{\infty}q^{2n^2+\frac{3n}2}
\frac{(q,q^{\frac12},q^{\frac12};q)_n(q,q^{\frac32},q^{\frac32};q)_{2n}(1-q^{2n+1})(1-q^{4n+2})}
     {(q^{\frac32};q)_{n}(q^{\frac32};q)_{2n}(q,q^{\frac32};q)_{3n}(1-q^{3n+1})(1-q^{3n+\frac32})^2}\\
&\times\bigg\{q^{n+\frac12}+
\frac{(1-q^{n+\frac12}) (1-q^{3n+1})(1-q^{3n+\frac32})^2 (1-q^{5 n+1})}
     {(1-q^{2 n+1})(1-q^{2n+\frac{1}{2}})^2(1-q^{3n+\frac{1}{2}})(1-q^{4n+2})}\\
&\qquad+q^{3n+2}
\frac{(1-q^{n+\frac{1}{2}})^2 (1-q^{n+1})(1-q^{3 n+\frac{3}{2}})(1-q^{5 n+4})}
     {(1-q^{3n+2})(1-q^{2 n+\frac32})(1-q^{3n+\frac{5}{2}})^2(1-q^{4n+2})}\bigg\}.}
\end{exam}
This provides a $q$-analogue for another series:
\pnq{\frac{9\pi^2}{16}
&=\sum_{n=0}^{\infty}
\dfrac{(2n)!(4n+1)!}
     {(6n+1)!(1+2n)(1+3n)}\\
&\times\bigg\{1+\frac{9 (1+2 n) (1+3 n) (1+5 n)}{2 (1+4 n)^2 (1+6 n)}
+\frac{3(1+n) (1+2 n)^2 (4+5 n)}{2 (2+3 n) (3+4 n) (5+6 n)^2}\bigg\}.}




\end{document}